\documentclass[hidelinks,onefignum,onetabnum]{siamart220329}

\usepackage[letterpaper]{geometry}

\usepackage{hyperref}
\usepackage{amsmath,amssymb,amsfonts}
\usepackage{eqnarray}
\usepackage{xcolor}
\usepackage{caption,subcaption,subfloat}
\usepackage{graphicx,float}
\usepackage[noend]{algpseudocode}

\usepackage{epstopdf}

\ifpdf
\DeclareGraphicsExtensions{.eps,.pdf,.png,.jpg}
\else
\DeclareGraphicsExtensions{.eps}
\fi

\newcommand{\eqt}[1]{Equation~(\ref{#1})}
\newcommand{\fig}[1]{Fig.~(\ref{#1})}
\newcommand{\sect}[1]{Section~(\ref{#1})}
\newcommand{\algo}[1]{Algorithm~(\ref{#1})}
\newcommand{\floor}[1]{\left \lfloor #1 \right \rfloor} 

\usepackage{wrapfig}
\usepackage{pgfplots}

\usepackage{tikz}
\usetikzlibrary{arrows}
\usetikzlibrary{decorations.shapes}
\usetikzlibrary{arrows.meta}
\usetikzlibrary{calc}

\tikzset{
	myarrow/.style={-{Triangle[length=3mm,width=1mm]}}
}

\usepackage[norndcorners,customcolors,nofill]{hf-tikz}
\hfsetbordercolor{black!50}

\DeclareMathOperator*{\argmin}{arg\,min}

\newcommand{\dimension}[1]{#1D}

\begin{document}

\newcommand\relatedversion{}
\renewcommand\relatedversion{\thanks{The full version of the paper can be accessed at \protect\url{https://arxiv.org/abs/1902.09310}}} 


\headers{Scalable Functional Approximation of Discrete Data}{V.S. Mahadevan, D. Lenz, I. Grindeanu and T. Peterka}

\title{Parallel Domain Decomposition Techniques Applied to Multivariate Functional Approximation of Discrete Data
\thanks{Submitted to editors \date{}.
\funding{This work was funded by the Early Career Research Program, Department of Energy, USA under contract no.~DE-AC02-06CH11357.}}
}

\author{Vijay S. Mahadevan\thanks{Argonne National Laboratory, Lemont, IL (\email{mahadevan@anl.gov})} 
\and
David Lenz\thanks{Argonne National Laboratory, Lemont, IL (\email{dlenz@anl.gov})}
\and
Iulian Grindeanu\thanks{Argonne National Laboratory, Lemont, IL (\email{iulian@anl.gov})}
\and
Thomas Peterka\thanks{Argonne National Laboratory, Lemont, IL  (\email{tpeterka@mcs.anl.gov})}
}

\date{}

\maketitle


\begin{abstract} \small\baselineskip=5pt 
	Compactly expressing large-scale datasets through Multivariate Functional Approximations (MFA) can be critically important for analysis and visualization to drive scientific discovery. Tackling such problems requires scalable data partitioning approaches to compute MFA representations in amenable wall clock times. We introduce a fully parallel scheme to reduce the total work per task in combination with an overlapping additive Schwarz-based iterative scheme to compute MFA with a tensor expansion of B-spline bases, while preserving full degree continuity across subdomain boundaries. While previous work on MFA has been successfully proven to be effective, the computational complexity of encoding large datasets on a single process can be severely prohibitive. Parallel algorithms for generating reconstructions from the MFA have had to rely on post-processing techniques to blend discontinuities across subdomain boundaries. In contrast, a robust constrained minimization infrastructure to impose higher-order continuity directly on the MFA representation is presented here. We demonstrate the effectiveness of the parallel approach with domain decomposition solvers, to minimize the subdomain error residuals of the decoded MFA, and more specifically to recover continuity across non-matching boundaries at scale. The analysis of the presented scheme for analytical and scientific datasets in 1-, 2- and 3-dimensions are presented. Extensive strong and weak scalability performances are also demonstrated for large-scale datasets to evaluate the parallel speedup of the MPI-based algorithm implementation on leadership computing machines.
\end{abstract}





\section{Introduction}
\label{sec:introduction}

Large-scale discrete data analysis of various scientific computational simulations often requires high-order continuous functional representations that have to be evaluated anywhere in the domain. Such expansions described as {\em Multivariate Functional Approximations} (MFA) \cite{de1983approximation, functional-analysis} in arbitrary dimensions allow the original discrete data to be compressed, and expressed in a compact form, in addition to supporting higher-order derivative queries (without further approximations such as finite differences) for complex data analysis tasks. MFA utilizes approximations of the raw discrete data using a hypervolume of piecewise continuous functions. One particular option is to use the variations of the B-Spline or NURBS bases \cite{nurbs-book, peterka-mfa} for the MFA {\em encoding} of scientific data. The reconstructed data in MFA retains the spatiotemporal contiguity, and statistical distributions, with lesser storage requirements. Due to the potentially large datasets that need to be encoded into MFA, the need for computationally efficient algorithms (in both time and memory) to parallelize the work is critically important. It is also essential to guarantee that the solution smoothness in the reconstructed (or {\em decoded}) dataset is consistently preserved when transitioning from a single MFA domain to multiple domains during parallelization.

Achieving improved performance without sacrificing discretization accuracy requires an infrastructure that is consistent in the error metrics of the decoded data and an algorithm that remains efficient in the limit of large number of parallel tasks. In this paper, we will utilize domain decomposition (DD) techniques \cite{smith-ddm} with data partitioning strategies to produce scalable MFA computation algorithms tha minimizes the reconstruction error when reproducing a given dataset. In such partitioned analysis, it is imperative to ensure that the continuity of the encoded and decoded data across subdomain interfaces is maintained, and remain consistent with the degree of underlying expansion bases used in MFA \cite{peterka-mfa}.  This is due to the fact that independently computing MFA approximations in individual subdomains do not guarantee even $C^0$ regularity in either the MFA space or in the reconstructed data. 
In order to tackle this issue, we rely on an iterative Schwarz-type DD scheme to ensure that continuity is enforced, and the overall error stays bounded as the number of subdomains are increased (or as the subdomain size decreases).

In addition to remaining efficient, we also require the devised algorithms to extend naturally to arbitrary dimensional settings and to handle large datasets. We next discuss some of the related work in the literature that have been explored for reconstruction of scattered data, and approaches to make these algorithms scalable in order to motivate the ideas presented in the paper. 




\subsection*{Literature Review}
\label{sec:related-work}

Domain decomposition (DD) techniques in general rely on the idea of splitting a larger domain of interest into smaller partitions or subdomains, which results in coupled Degrees-of-Freedom (DoF) at their common interfaces. Typical applications of DD in Boundary-Value problems (BVP) \cite{smith-ddm, lions-asm} have been successfully employed to efficiently compute the solution of large, discretized Partial Differential Equations (PDEs) in a scalable manner. 
DD techniques for parallel approximation of scattered data have been explored previously with Radial Basis Functions (RBF) \cite{mai-approx-rbf}, yielding good scalability and closely recovering the underlying solution profiles. In general, overlapping multiplicative and additive Schwarz \cite{orasm-as-ms-2007} iterative techniques for RBF \cite{ddm-rbf} have proven successful to tackle large-scale problems. Additionally, the use of restricted variants of additive Schwarz (RAS) method as preconditioners, with Krylov iterative solvers, can yield iterative schemes \cite{yokota-rasm-rbf} with $O(N)$ computational complexity, as opposed to the typical $O(N log(N))$ complexity with traditional RBF reconstructions \cite{ddm-rbf-fast}. The extensions of these ideas to B-spline bases exposes a way to fully parallelize traditional, serial MFA computations.

Combining the application of DD schemes and NURBS bases with isogeometric analysis (IGA) \cite{cottrell2009, da2012} for high-fidelity modeling of nonlinear PDEs \cite{dede2015, marini2015parallel, petiga-dalcin-2016} have enjoyed recent success at scale. However, many of these implementations lack full support to handle multiple geometric patches in a distributed memory setting due to non-trivial requirements on continuity constraints at patch boundaries. 
Directly imposing higher-order geometric continuity in IGA requires specialized parameterizations in order to preserve the approximation properties \cite{kapl2018construction}, which can be difficult to parallelize \cite{hofer2018fast} generally. 
In a similar vein, using B-spline bases to compute the MFA in parallel, while maintaining higher-order continuity across subdomains has not been fully explored previously. 

To overcome some of these issues with discontinuities along NURBS or B-spline patches, Zhang et al. \cite{zhang-nurbs-continuity} proposed to use a gradient projection scheme to constrain the value ($C^0$), the gradient ($C^1$), and the Hessian ($C^2$) at a small number of test points for optimal shape recovery. Such a constrained projection yields coupled systems of equations for control point data for local patches, and results in a global minimization problem that needs to be solved.

Alternatively, it is possible to create a constrained recovery during the actual post-processing stage i.e., during the decoding stage of the MFA through standard blending techniques \cite{grindeanu-blending}, in order to recover continuity in the decoded data. However, the underlying MFA representation remains discontinuous, and would become more so with increasing number of subdomains without the ability to recover higher-order derivatives along these boundaries. Moreover, selecting the amount of overlaps and resulting width of the blending region relies strongly on a heuristic, which can be problematic for general problem settings.

In contrast, we propose extensions to the constrained solvers used by Zhang et al. \cite{zhang-nurbs-continuity} and Xu et al. \cite{xu-jahn-discrete-adjoint}, and introduce a two-level, DD-based, parallel iterative scheme to enforce the true degree of continuity, independent of the basis function polynomial degree $p$, unlike the low-order constraints used previously  \cite{zhang-nurbs-continuity}. The outer iteration utilizes the RAS method \cite{gander-rasm}, with efficient inner subdomain solvers that can handle linear Least-Squares systems to minimize the decoded residual within acceptable error tolerances. 
Such an iterative solver has low memory requirements that scales weakly with growing number of subdomains, and necessitates only nearest-neighbor communication of the interface data once per outer iteration to converge towards consistent MFA solutions.

\subsection*{Structure of the paper}

The paper is organized as follows. \sect{sec:approach} presents the theory and necessary details about the subdomain solvers, and the DD approach used to converge the boundary continuities across MFA subdomains. Next, in \sect{sec:results}, the DD solver is applied to several \dimension{1}, \dimension{2} and \dimension{3} synthetic and real-world datasets to verify error convergence, and the parallel scalability of the iterative algorithm for decreasing subdomain sizes is demonstrated. Finally, key observations from the parallel MFA solver and future extensions to more complex cases with spatial adaptivity are presented in \sect{sec:conclusions}.



\section{Approach}
\label{sec:approach}

With motivations to accelerate the computation of an accurate MFA representation scalably, we utilize a data decomposition approach with overlapping subdomains to create shared layers of piecewise accurate functional reconstructions. This is similar to a multipatch approach typically taken in IGA computations \cite{cottrell2009, petiga-dalcin-2016}.  However, in order to ensure that higher-order continuity across domain boundaries are preserved, an outer iteration loop is inevitable to converge the shared unknowns across the interfaces. These global iterations guarantee consistent MFA encodings in parallel, without which the representations will not even ensure $C^0$ regularity. 

In this section, we first provide an illustrative example by formulating the constrained minimization problem to be solved in each subdomain and explain the iterative methodology used in the current work to converge the shared DoFs. We will also introduce the idea of using open vs closed knots, which are clamped or floating respectively at subdomain boundaries and discuss the advantages of using one approach over the other. 


\subsection{Numerical Background}
\label{sec:background}

A $p$-th degree NURBS or B-spline curve \cite{nurbs-book} is defined using the Cox-deBoor functions for each subdomain as

\begin{eqnarray}
	\vec{C}(u) &=& \sum_{i=0}^{n} R_{i,p}(u) \vec{P}(i), \quad \forall u \in \Omega \\
	R_{i,p}(u) &=& \frac{N_{i,p}(u) W_i}{\sum_{i=0}^{n} N_{i,p}(u) W_i}
	\label{eq:nurbs-basis}
\end{eqnarray}

where $R_{i,p}(u)$ are the piecewise rational functions with $\vec{P}$ control points of size $n$, $W_i$ are the control point weights, with the $p$-th degree B-spline bases $N_{i,p}(u)$ defined on a knot-vector $u$. Note that exact high-order derivatives of these B-spline basis defined in \eqt{eq:nurbs-basis} can also be evaluated without any approximation errors at the control point locations using the Cox-deBoor recurrence relations \cite{de1983approximation}. This property becomes especially important when performing analysis and in-situ visualization directly based on the MFA representation of underlying data \cite{mfa-vis-dvr}.

Given a set of input points $\vec{Q}$ that need to be encoded into a MFA, with the weights $W=1$ (B-spline representations) for simplicity, the unconstrained minimization problem to compute the optimal set of control point locations within a subdomain can be posed as a solution to a linear Least-SQuares (LSQ) system that minimizes the net error of the B-spline approximation.

\begin{eqnarray}
	\argmin_{\vec{P} \in \mathbb{R}^n} {E} = \left\lVert \vec{Q} - R \vec{P} \right\rVert_{L_2}, \quad \quad R \in \mathbb{R}^{m \times n}, \vec{Q} \in \mathbb{R}^m
	\label{eq:minimization-problem}
\end{eqnarray}

An appropriate LSQ solver such as the one based on Cholesky decomposition or the more efficient $\ell$-BFGS scheme \cite{zheng-bo-bspline-bfgs} can compute the control point solution $\vec{P}$ that minimizes the residual error $\vec{E}$ for the given input data $\vec{Q}$ and MFA representation of degree $p$. Note that the minimization procedure can be performed independently on each subdomain without dependencies as there are no constraints explicitly specified in \eqt{eq:minimization-problem}.
However, in order to recover high-order continuity across subdomain interfaces, computing unconstrained solutions is insufficient. At a minimum, the DoFs lying on the shared subdomain boundaries have to be converged to recover $C^0$ continuity for the decoded solution data ($R \vec{P}$).
%

More generally, the constrained minimization problem to recover continuity \cite{nurbs-book} can be formulated as 
\begin{equation}
	R \vec{P} = \vec{Q} \quad \mid \quad \mathcal{C} \vec{P} = \vec{G}, \label{eq:global-constrained-problem}
\end{equation}
where $\mathcal{C}$ is the constraint matrix imposing continuity restrictions on the control points $\vec{P}$ along with its derivatives, with data exchanged from neighboring domains stored in $\vec{G}$, around the neighborhood of the interface $\Omega_{i,j}$ shared by subdomains $i$ and $j$.
With the use of penalized constraints ($\mathcal{C}$) and Lagrange multipliers \cite{dornisch2011, paul2020}, the solution to the constrained LSQ problem can recover optimal control point values. 

A straightforward approach to achieve $C^0$ continuity in the recovered solution is by ensuring that the common control point data $\vec{P}$ at subdomain interfaces are clamped with repeated knots, in addition to using clamping at the global domain boundaries. In this scheme, the control points exactly interpolate (are clamped to) input data points at the subdomain interface boundaries. Such an approach requires in general a good spatial distribution of $\vec{Q}$, and yields only low-order continuous approximations ($C^0$) when the solution remains smooth across the subdomain interfaces. It should also be noted that as the number of subdomains increases, the global solution being computed becomes further constrained, and more interpolatory due to clamped DoFs. Moreover, the MFA solution computed becomes dependent on the number of subdomains used to decompose the problem; i.e., the global control point data $\vec{P}$ recovers different reconstructions as a function of number of subdomains ($\mathcal{N}$) used.


While the numerics and implementation of the domain decomposed MFA can be much simpler with clamped knots on all subdomain boundaries, ensuring higher-order continuity would require that all $p-1$ derivatives of the approximation match as well. As a continuous extension, one could relax the interpolatory behavior of clamped knot boundaries by reducing the number of repeated knots, and instead use floating (or unclamped) knots at internal subdomain boundary interfaces by sharing knot spans between subdomains. This modification allows us to recover fully consistent ($C^0$ to $C^{p-1}$) continuous MFA reconstructions using the solution procedure detailed for the global constrained minimization problem \eqt{eq:global-constrained-problem}.

\subsection{Shared Knot Spans at Subdomain Interfaces}

\begin{figure}[htbp]
	\centering
	\subfloat[Even degree $p=2$\label{fig:degree-2-1d}]{%
		\includegraphics[width=0.47\textwidth]{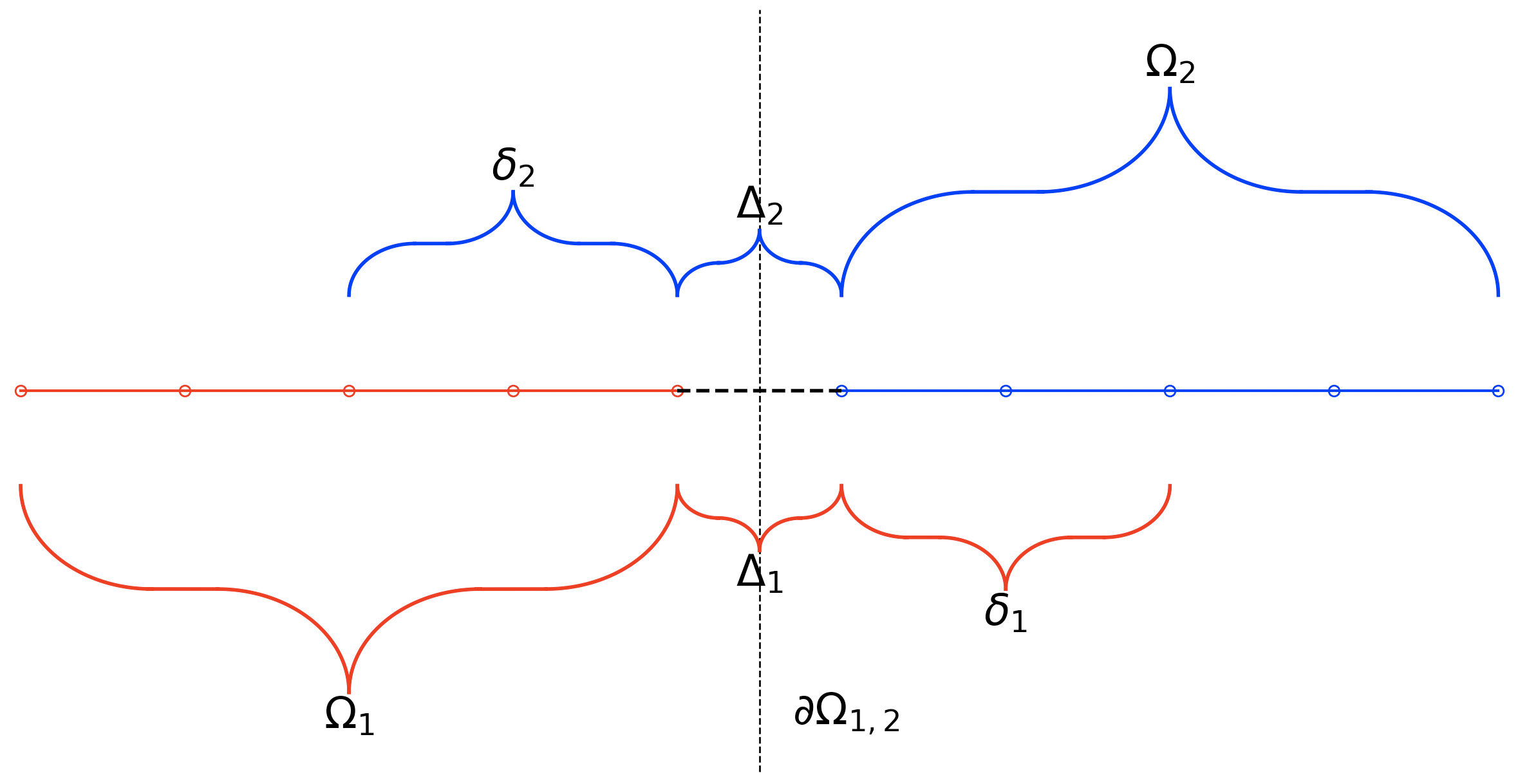}}
	\hfill
	\subfloat[Odd degree $p=3$\label{fig:degree-3-1d}]{%
		\includegraphics[width=0.47\textwidth]{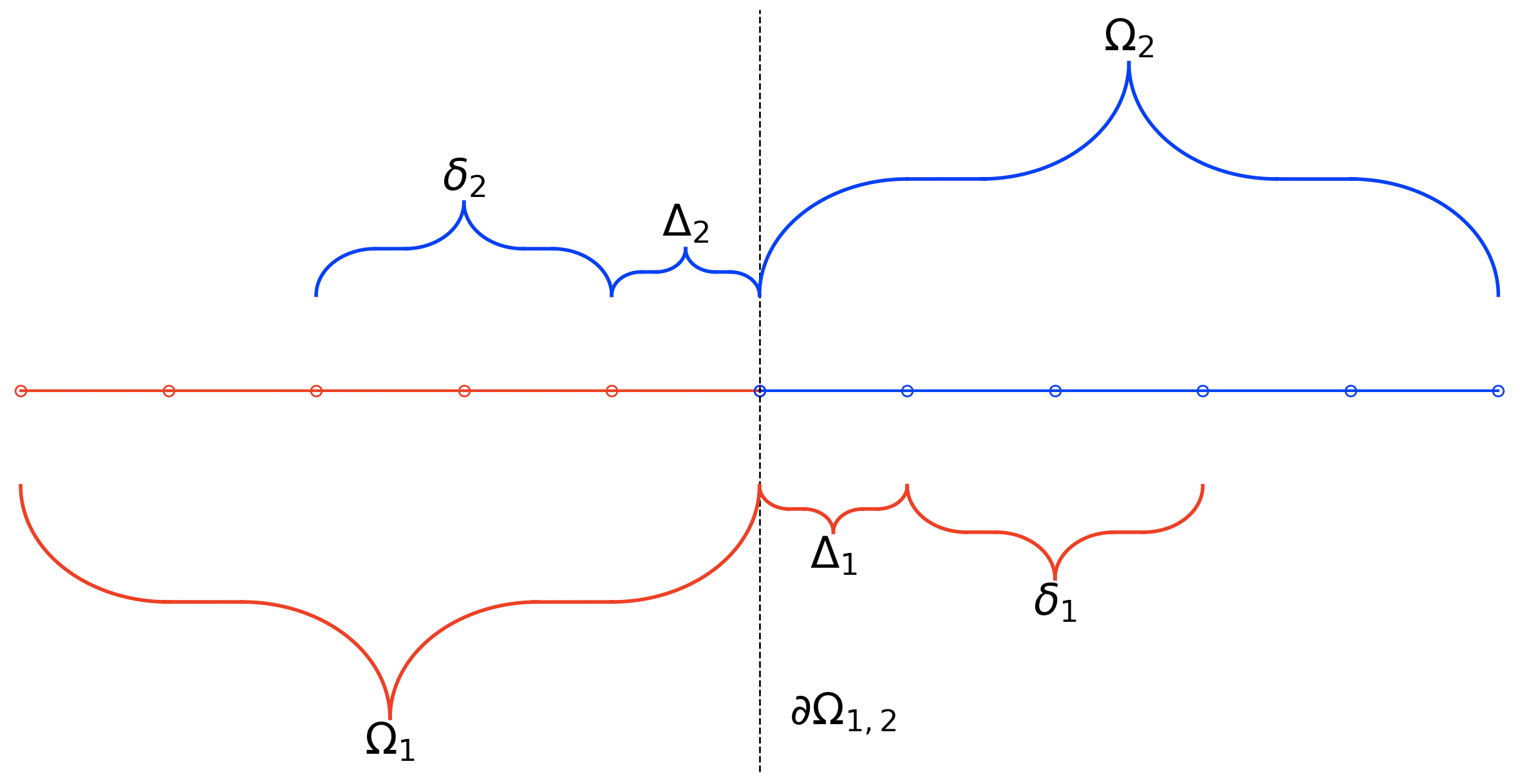}}
	\caption{Illustration: \dimension{1} parallel partitioned domain with floating (unclamped) interior knots and augmented spans ($\left| \delta \right|=2$)}
	\label{fig:DD-subdomain-illustration}
\end{figure}


Instead of using clamped knots, we utilize floating (unclamped), shared knot spans near all interior subdomains such that the high-order continuity and consistency of the reconstructed solution with respect to $\mathcal{N}$ are preserved. 

For the purpose of illustration and to explain the proposed solver methodology, let us consider a simple one dimensional domain ($\Omega$) with two subdomains ($\mathcal{N}=2$) as shown in \fig{fig:DD-subdomain-illustration}, where $\Omega_1$ and $\Omega_2$ represent the subdomains that share an interface $\partial \Omega_{1,2}$. In \fig{fig:DD-subdomain-illustration}, the layout of the knot spans for both an even degree ($p=2$) and odd degree ($p=3$) are shown. For generality, we also introduce here an overlap layer $\Delta_1$ and $\Delta_2$ on each subdomain that represents the set of shared knot spans with its adjacent subdomain (for internal boundaries), and an optional augmented layer $\delta_1$ and $\delta_2$ that has a connotation similar to that of an overlap region in traditional DD schemes \cite{smith-ddm}. Note that in order to reconstruct the input data in $\Omega_i, \forall i \in [1,2]$, the knot spans must mandatorily include $\Delta_i$ regions. This $\Delta_i$ overlap region is required by definition to maintain partition of unity of a B-spline curve in order to evaluate \eqt{eq:nurbs-basis}. For generality, $\Delta_i$ represents the repeated knots along clamped global domain boundaries, and the shared knots between two subdomains in the unclamped interior boundaries.
For arbitrary degree $p$, the number of knot spans in $\Delta_i$ is given by $\floor{\frac{p}{2}}$, where $\floor{.}$ represents the floor operator. In multidimensional tensor product expansions, these shared spans are replaced by shared layers of knot spans along the subdomain interfaces. The $\delta_i$ regions are additional, and optional, shared knot spans that can help improve error convergence in a manner similar to overlap regions in DD methods used for PDE solvers \cite{gander-rasm, smith-ddm}. 

The control point DoF vector can be represented by three separate parts based on the local support of the basis expansion. The control point vector is in general given as $\vec{P} = [\vec{P}(\Omega); \vec{P}(\Delta); \vec{P}(\delta)]$. For the \dimension{1} scenario illustrated, this is shown below for $\mathcal{N}=2$, $p=3$ and $\left| \delta \right|=1$, where the operator $\left| \mathcal{D} \right|$ represents the number of knot spans in any underlying domain $\mathcal{D}$. 
\vspace*{-1mm}
\begin{equation}
	\tikzset{offset def/.style={
			above left offset={-0.1,0.8},
			below right offset={0.1,-0.65},
		}
	}
	\hfsetfillcolor{red!10}
	\hfsetbordercolor{red}
	\vec{P}_1 =
	\left[
	\begin{array}{c}
		\tikzmarkin{A}(1,-0.2)(-0.1,0.3) P_1(1)   \\
		P_1 (2) \\
		\vdots   \\
		\tikzmarkend{A} P_1 (m) \\
		$\quad$ \vspace*{-2mm} \\
		\tikzmarkin{B}(0.1,-0.2)(-0.1,0.3) P_2(1) \tikzmarkend{B} \\
		$\quad$ \vspace*{-2mm} \\
		\tikzmarkin{Bp}(0.1,-0.2)(-0.1,0.3) P_2(2) \tikzmarkend{Bp} \\
	\end{array}
	\right]
	, \quad \quad
	\hfsetfillcolor{blue!10}
	\hfsetbordercolor{blue}
	\vec{P}_2 =
	\left[
	\begin{array}{c}
		\tikzmarkin{D}(0.1,-0.2)(-0.1,0.3) P_2(1)   \\
		P_2 (2) \\
		\vdots   \\
		P_2 (n) \tikzmarkend{D} \\
		$\quad$ \vspace*{-2mm} \\
		\tikzmarkin{Cp}(0.1,-0.2)(-0.1,0.3) P_1(m) \tikzmarkend{Cp} \\
		$\quad$ \vspace*{-2mm} \\
		\tikzmarkin{C}(0.1,-0.2)(-0.1,0.3) P_1(m-1) \tikzmarkend{C} \\
	\end{array}
	\right]
	\label{eqn:vecp-layout}
\end{equation}

\begin{tikzpicture}[remember picture,overlay]
	\pgfsetarrowsend{latex} 
	%
	\coordinate (A-aa) at ($(A)+(-0.8,-0.7)$);
	\node[align=left,left] at (A-aa) {\footnotesize{$\vec{P}_1(\Omega_1)$}};
	\path[>=stealth,red,draw] (A-aa) -- ($(A)+(0.1,-0.7)$);
	%
	\coordinate (B-aa) at ($(B)+(-0.8,-0.2)$);
	\node[align=left,left] at (B-aa) {\footnotesize{$\vec{P}_1(\Delta_1)$}};
	\path[>=stealth,blue,draw] (B-aa) -- ($(B)+(0.1,-0.2)$);
	%
	\coordinate (Bp-aa) at ($(Bp)+(-0.8,-0.2)$);
	\node[align=left,left] at (Bp-aa) {\footnotesize{$\vec{P}_1(\delta_1)$}};
	\path[>=stealth,blue,draw] (Bp-aa) -- ($(Bp)+(0.1,-0.2)$);
	%
	\coordinate (D-aa) at ($(D)+(2.5,-1.0)$);
	\node[align=right,right] at (D-aa) {\footnotesize{$\vec{P}_2(\Omega_2)$}};
	\path[>=stealth,blue,draw] (D-aa) -- ($(D)+(0.9,-1.0)$);
	%
	\coordinate (Cp-aa) at ($(Cp)+(2.5,-0.3)$);
	\node[align=right,right] at (Cp-aa) {\footnotesize{$\vec{P}_2(\Delta_2)$}};
	\path[>=stealth,red,draw] (Cp-aa) -- ($(Cp)+(1.0,-0.3)$);
	%
	\coordinate (C-aa) at ($(C)+(2.8,-0.3)$);
	\node[align=right,right] at (C-aa) {\footnotesize{$\vec{P}_2(\delta_2)$}};
	\path[>=stealth,red,draw] (C-aa) -- ($(C)+(1.6,-0.3)$);
\end{tikzpicture}

where $m, n$ are the number of control points in $\Omega_1$ and $\Omega_2$ respectively. Note that higher degree expansions (for e.g., $p>3$) will require more support points in $\Delta$ from adjacent subdomains in order to decode the MFA up to $\partial \Omega_{1,2}$. This implies that $P_i(\Delta_i)$ in addition to the optional $P_i(\delta_i)$ vectors directly provide a measure of the required cost of communication with adjacent subdomains.

Now, the constrained minimization problem for the two subdomain case can be written as
\begin{equation}
	\left[
	\begin{array}{c|c}
		R_{1}(\Omega_1) & \lambda_{1,2}(\Delta_1 \cup \delta_1) \\
		\hline
		\lambda_{2,1}(\Delta_2 \cup \delta_2) & R_{2}(\Omega_2)
	\end{array}
	\right]
	\left[
	\begin{array}{c}
		\vec{P}_{1} \\
		\vec{P}_{2}
	\end{array}
	\right]
	=
	\left[
	\begin{array}{c}
		\vec{Q}_{1} \\
		\vec{Q}_{2}
	\end{array}
	\right]
	\label{eq:global-system}
\end{equation}


where the diagonal operators $R_{1}$ and $R_{2}$ are the piecewise rational functions that minimize the local subdomain residuals in $\Omega_j, \forall j \in [1,2]$, while the off-diagonal blocks $\lambda_{1,2}$ and $\lambda_{2,1}$ represent the coupling terms between the subdomains near the interface $\partial \Omega_{1,2}$. This coupling term provides the constraints on the shared control point data, and higher-order derivatives as needed to recover smoothness and enforce continuity along subdomain boundaries. For higher dimensional problems, the constraints on the control points must include both face neighbor and diagonal neighbor contributions to accurately determine the globally consistent minimization problem. 

The coupling blocks $\lambda_{i,j}$ can be viewed as Lagrange multipliers that explicitly couple the control point DoFs across a subdomain interface ($\vec{P}_{1} \cap \vec{P}_{2}$) such that continuity is preserved in a weak sense \cite{nurbs-book}. Using appropriate Schur complements to eliminate the coupled DoF contributions in each subdomain, with $\lambda_{i,j}$ evaluated at \textit{lagged} iterates of adjacent subdomains, the set of coupled constrained equations in \eqt{eq:global-system} can be completely decoupled for each subdomain. This modified system resembles a block-Jacobi operator of the global system. The scheme illustrated in this section follows ideas similar to the Jacobi-Schwarz method \cite{gander-rasm} and the overlapping,  restricted-Additive-Schwarz (RAS) scheme \cite{orasm-as-ms-2007}.

In the above description, the coupled data chunks, $\vec{P}_1(\Delta_1)$ and $\vec{P}_2(\Delta_2)$ belonging to adjacent subdomains near $\partial \Omega_{1,2}$ are exchanged simultaneously before the local domain solves are computed. 
One key advantage with such a DD scheme is that only nearest neighbor exchange of data is required, which keeps communication costs bounded as the number of subdomains $\mathcal{N}$ increase \cite{orasm-as-ms-2007, gander-rasm}, while providing opportunities to interlace recomputation of the constrained control point solution. Note that in a RAS iterative scheme, nearest neighbor exchanges can be performed compactly per dimension and direction, thereby minimizing communication costs and eliminating expensive global collectives.

\subsubsection*{Augmenting Knot Spans with Overlap}

One of the key metrics of interest is that the parallel solver infrastructure does not amplify any approximation errors unresolved by the tensor product B-spline mesh. Since the local decoupled subdomain solution is encoded accurately to satisfy \eqt{eq:minimization-problem} in each individual subdomain without any data communication (i.e., embarassingly parallel), imposing the constraints for the shared DoFs in $\Delta$ should ensure the error change is bounded. However, as the control point data across subdomains become synchronized, numerical artifacts, especially for high-degree ($p>2$) basis reconstructions at subdomain interfaces can become dominant sources of error. A key metric to consider in all experiments is to validate that the multiple subdomain case produces the same error profile as a single subdomain case, in order to ensure convergence of the solvers to the same unique solution, independent of $\mathcal{N}$. 

For many problem domains, overlapping variants of Schwarz solvers \cite{lions-asm, gander-rasm} have been proven to be more stable, efficient and scalable compared to nonoverlapping variants \cite{bjorstad-overlap-1989, orasm-as-ms-2007}. We utilize the concept of overlap regions by sharing additional knot spans between subdomains in order to produce better MFA reconstructions of the underlying data. This user-specified, additional overlap is described by $\delta_j, \forall j \in [1,2]$ in \fig{fig:DD-subdomain-illustration}. The amount of data overlap utilized for computing the functional approximation can directly affect the conditioning in the subdomain solver, and the scalability of the overall algorithm. Additionally, we expect the residual errors $\vec{E}$ from the approximation to remain bounded as the number of subdomain increase with appropriate overlap regions. 

\begin{figure}[htbp]
	\centering
	\begin{tikzpicture}
		\begin{scope}[xshift=1.5cm]
			\node[anchor=south west,inner sep=0] (image) at (0,0) {\includegraphics[width=0.44\textwidth]{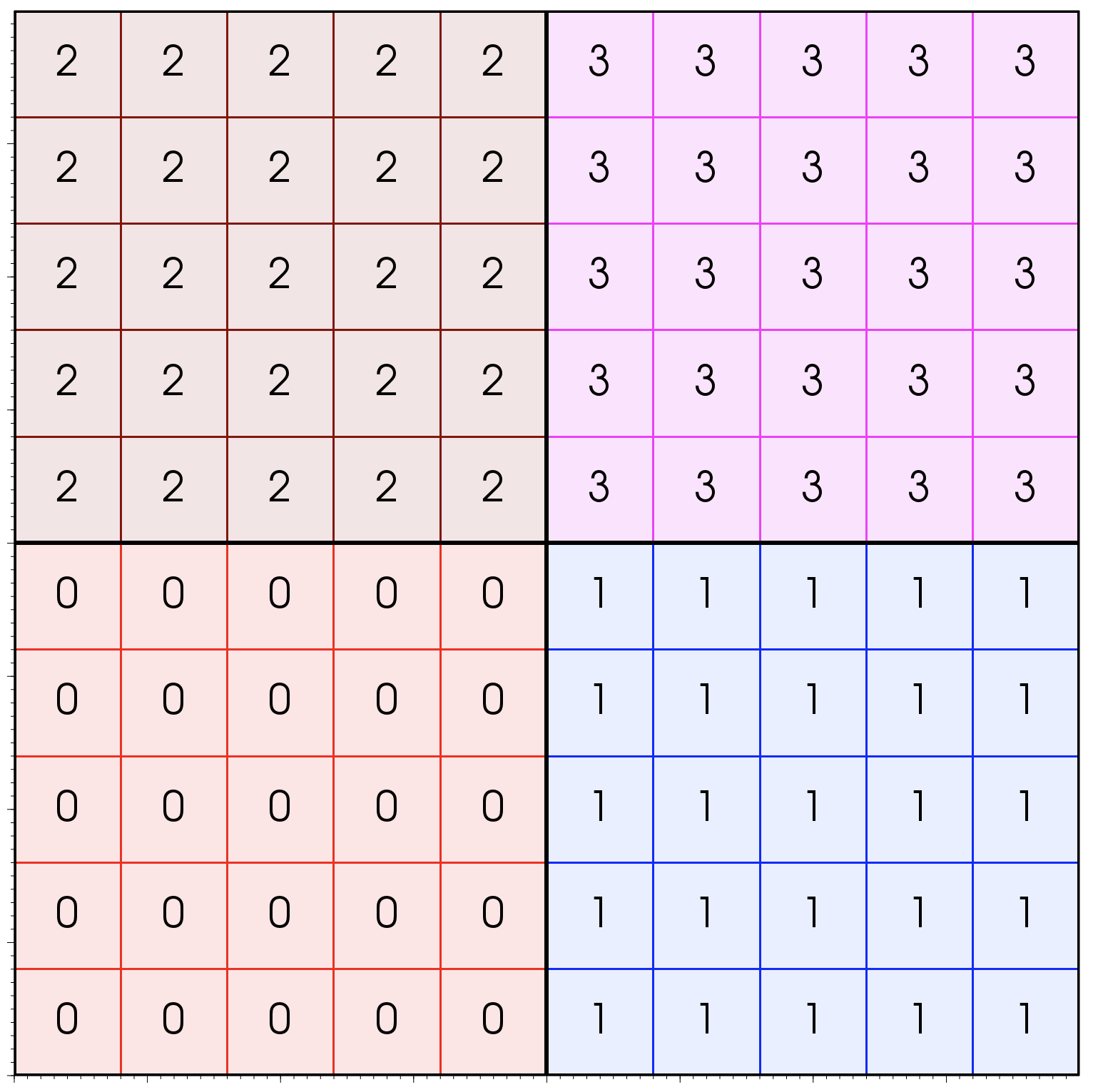}};
			\begin{scope}[x={(image.south east)},y={(image.north west)}]
				\draw[black,ultra thick,rounded corners] (0.49,0.0) rectangle (0.61,0.61);
				\draw[black,ultra thick,rounded corners] (0.0,0.49) rectangle (0.61,0.61);
				\draw[black,ultra thick,densely dashed,rounded corners] (0.59,0.0) rectangle (0.71,0.71);
				\draw[black,ultra thick,densely dashed,rounded corners] (0.0,0.59) rectangle (0.71,0.71);
				\draw[black,ultra thick,dotted,rounded corners] (0.0,0.0) rectangle (0.51,0.51);
				
				\draw[red,ultra thick,rounded corners] (0.29,0.29) rectangle (0.71,0.71);
				\node[draw] at (-0.07, 0.25) {\large $\Omega_0$};
				\node[draw] at (-0.07, 0.55) {\large $\Delta_0$};
				\node[draw] at (-0.07, 0.66) {\large $\delta_0$};
				\node[draw] at (0.55,-0.07) {\large $\Delta_0$};
				\node[draw] at (0.66,-0.07) {\large $\delta_0$};
				
				\draw [-stealth, line width=2pt, red] (0.71,0.5) -- (1.1,0.5);
				\node[draw] at (1.35,0.5) {multishared region};
			\end{scope}
			
		\end{scope}
	\end{tikzpicture}%
	\caption{\dimension{2} subdomains with $\mathcal{N}=4$, $p=3$ and the augmented overlap $\left| \delta \right|=1$ showing local subdomains $\Omega_i$, mandatory overlap for floating knots $\Delta_i$, optional overlap regions $\delta_i$, where $i \in [0,\mathcal{N}-1]$, and multishared DoF regions that couple multiple subdomains (marked in red).}
	\label{fig:2d-schematic}
\end{figure}
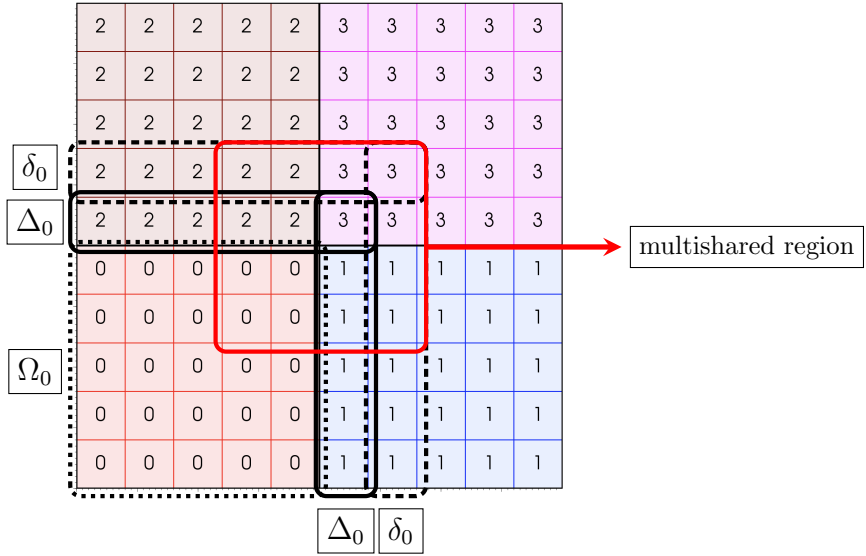

For better clarity, we will use overlap regions $\delta$, as illustrated in \fig{fig:2d-schematic} for a \dimension{2} problem with $p=3$ and $\mathcal{N}=4$, to increase the size of the local problem ($\Omega$), and to improve the accuracy of domain decomposed approximations. We note that the control point data $\vec{P}$ in both the $\Delta$ and $\delta$ overlap regions are shared and uniformly weighted (averaged) by $\vec{P}$ computed in the neighboring subdomains.

Note that in the \dimension{2} schematic, shared data in $\delta$ regions are always exchanged between neighboring subdomains. The set of $\vec{P}(\delta)$ are explicitly used only to impose constraints that contribute to the reconstruction of datasets, and hence play a role in the approximation error of MFA. It is also important to note that when $\vec{P}$ DoFs are {\em multishared} between subdomains, then shared data between multiple $\Omega_j, \forall j \in [0,3]$ need to be exchanged in order to compute $\vec{P}_i(\delta_i), \forall i \in [0,3]$.


%
%
%
%

We also introduce a definition for compression ratio ($\eta$), which gives the ratio of the total input points in the dataset ($dim\text{ }\vec{Q}$) to the total control points ($dim\text{ }\vec{P}$) used in the MFA B-spline representation. As $\eta \rightarrow 1$, one can achieve better error residuals in comparison to the reference data for a given degree $p$, while $\eta \gg 1$ produces smooth approximations with larger error profiles. When we analyze the results in this manuscript, it should be remembered that the reconstruction error is always inversely proportional to $\eta$.

Next, the parallel MFA computation workflow that will be used with domain-decomposed subdomain partitions is presented in detail. 

\subsection{Solver Workflow}
\label{sec:solver-methodology}

Computing the functional approximation to large-scale datasets requires efficient solvers at two levels: firstly, the local decoupled subdomain solver for  \eqt{eq:minimization-problem}, and next, the constrained minimization problem in \eqt{eq:global-constrained-problem}. Hence, the global problem reduces to a series of local minimization problems in each subdomain.

\begin{equation}
	\begin{aligned}
		%
		& \argmin_{\vec{P} \in \mathbb{R}^n} & & \left\lVert \vec{Q}_{\ell} - R_{\ell} \vec{P}_{\ell} \right\rVert_{L_2}, \quad R_{\ell} \in \mathbb{R}^{m \times n}, \vec{Q}_{\ell} \in \mathbb{R}^m, \vec{P}_{\ell} \in \Omega_i \\
		& \text{subject to} & & \sum_{\partial \Omega_{i,j}} \left[ \mathcal{F}_{ij}( \vec{P}_{i}(\Delta_i), \vec{P}_{j}(\Omega_j) )  + \mathcal{F}_{ij}( \vec{P}_{i}(\delta_i), \vec{P}_{j}(\Omega_j) )  \right]^2 = 0, \quad \forall i, j \in [1,\ldots\mathcal{N}]
		%
		\label{eq:nonlinear-residuals}
	\end{aligned}
\end{equation}

where $\Omega_j$ are the neighboring subdomains of $\Omega_i$, $\mathcal{F}_{ij}(a,b)$ is the jump term across the shared interface DoFs $a$ and $b$ defined on subdomains $i$ and $j$ respectively. 

\subsubsection{Subdomain Solvers}

For the linear LSQ solvers that can be used to compute local subdomain control point solution $\vec{P}$, there are a variety of choices available. Direct methods like Singular Value decomposition or Cholesky decomposition operating on the normal equations \cite{bjorck1996} can compute optimal values. Alternatively, the iterative LSQ solvers such as orthogonal decomposition methods based on QR and QZ factorizations are more stable, especially when the normal form of the operator, $R^T R$, is ill-conditioned. 

\subsubsection{Restricted Additive-Schwarz Solvers}

The outer RAS iterations work together with nearest neighbor communication procedures to exchange shared DoF data between adjacent subdomains. This is an important step to ensure that $\vec{P}$ data computed through the LSQ procedure are consistent and high-order continuous across subdomain boundaries. The final minimized control point solution is achieved when the interface solutions match on all $\partial \Omega_{i,j} \in \Omega$ rendering zero jump residuals ($\mathcal{F}_{ij}$) on $\Delta$ and $\delta$ shared domains in \eqt{eq:nonlinear-residuals}.

It is also important to note that unlike the blending approaches that can be directly applied on decoded data \cite{grindeanu-blending}, the numerical error with the constrained iterative scheme is not bounded by the original partitioned, unconstrained least-squares solution; i,e., imposing subdomain boundary constraints can create artificial numerical peaks (non-monotonic) in reconstructed data as we converge towards continuity recovery. A solution to address this issue is to increase the amount of overlap range to ensure uniform convergence to the true single-subdomain solution error, even as the number of subdomains ($\mathcal{N}$) increases.


The nonoverlapping and overlapping RAS scheme applied to the computation of MFA exhibits scalable convergence properties in the limit of decreasing subdomain size (i.e., as $\mathcal{N} \to \infty$). This is a favorable property for strong scaling, especially when tackling large datasets, as the net computational cost always remains bounded. This behavior can be explained by the nature of how the RAS iterative procedure resolves the shared DoFs.

By using a weighted averaging procedure for all shared DoFs that reside in the $\Delta$ and $\delta$ domains, each outer iteration resolves any disparity in $\vec{P}$. The DoFs values on shared vertices ($d > 0$), edges ($d > 1$) and faces ($d > 2$) are resolved in the following order in consequent RAS iterations.
\begin{enumerate}
	\item Singly-shared ($\mathcal{SS}$) DoFs that are shared between two adjacent, neighboring subdomains; e.g., direct interface data belonging to $(\Delta_i \cup \delta_i) \cap \Omega_j$ in \fig{fig:2d-schematic}, $\forall i \in [1, \mathcal{N}]$ and $j \in [1, \hat{\mathcal{N}}]$, where $\hat{\mathcal{N}}$ represents the set of nearest neighbors sharing an interface $\Omega_{i,j}$
	\item Multi-shared ($\mathcal{MS}$) DoFs that are shared by multiple (more than 2) neighboring subdomains; e.g., diagonal corners that result from $\mathcal{S}_0 \cap \mathcal{S}_1 \cap \mathcal{S}_2 \cap \mathcal{S}_3$ in \fig{fig:2d-schematic}, where $\mathcal{S}_i=\Omega_i\cup\Delta_i\cup\delta_i$.
\end{enumerate}
Given both these specific DoF groups, the overlapping RAS scheme applied to MFA computation {\em always converges in 2 outer iterations} for problems in \dimension{2} and \dimension{3}, and a single iteration in \dimension{1} (due to lack of $\mathcal{MS}$ DoFs). This result is demonstrated in \sect{sec:results}. Note that achieving full convergence with high-order continuity in a single iteration is possible by combining constraints for both $\mathcal{SS}$ and $\mathcal{MS}$ DoFs. However, due to complexities in the implementation of constraint matching in \eqt{eq:nonlinear-residuals} for this algorithmic optimization, it has not been pursued here.

It is necessary to mention that we use a uniform weighting procedure to converge the shared DoFs between different subdomains, where the weights for each shared DoF is assigned as $w_i=\frac{1}{n_s}$, with $n_s$ being the number of subdomains containing DoF $i$ within its domain $\mathcal{S}$. This weighing procedure can be trivially replaced by Shephard's functions, especially in the context of adaptive discretizations with variable knot displacements.


\subsubsection{Note on Performance Characteristics}

%
%
The volume of messages exchanged between subdomains depends on several computational factors.

\begin{enumerate}
	\item \textbf{Clamping}: If the boundary knots are pinned, or if they have a floating knot description at subdomain interfaces depending on whether $C^0$ or $C^{p-1}$ continuity is required,
	\item \textbf{Parity}: Whether the MFA degree of expansion is odd or even, which determines the range of common knot spans shared between adjacent domains as given by $\Delta$ (refer to \fig{fig:DD-subdomain-illustration}),
	\item \textbf{Overlap}: The amount of augmented overlap ($\delta$), which determines the number of additional coupled data layers to be communicated between neighboring domains, both in terms of the input span space $\vec{Q}$, and control point DoFs $\vec{P}$ (refer to \fig{fig:DD-subdomain-illustration} and \fig{fig:2d-schematic}).
\end{enumerate}

At convergence, the interface data at $\partial \Omega_{1,2}$ will satisfy the higher-order continuity prescriptions specified by the user, thereby guaranteeing full regularity of $C^{p-1}$. The illustration in \fig{fig:DD-subdomain-illustration}, and the methodology description in this section can be generalized and extended to arbitrary dimensions in the tensor-product setting (with the parametric domain represented by a $d-$dimensional hypercube) as shown in \fig{fig:2d-schematic} for \dimension{2}. Using knot insertion and removal strategies based on deCasteljau subdivision procedures \cite{nurbs-book}, individual subdomains can also be adapted to resolve fast varying solutions and to reduce decoded error to be within user-specified tolerances \cite{nashed-rational}. While adaptivity has not been fully explored in the current work, enabling variable resolutions in different dimensions is a natural extension of the work that will still preserve high-order continuity.
The implementation of the presented approach with domain decomposition strategies, combined with overlapping RAS scheme yields a scalable scheme that will be demonstrated to be suitable for tackling large-scale data analysis problems.

\subsection{Implementation}
\label{sec:implementation}

%
	%
	%

\begin{algorithm}
	\caption{Domain Decomposed MFA Solver}
	\label{alg:pseudocode}
	\begin{algorithmic}
		\State{\textbf{Input:} Dataset and coordinates}
		\State{\textbf{Parameters:} $\mathcal{N}, p, m$}
		\State{\textbf{Setup:} Decompose domain into blocks with DIY, compute $R$}
		\While{$\vec{P}(\Omega)$ not converged}
		\For{$i \gets 1$ to $\mathcal{N}$}
		\If{$\vec{P}_i = 0$}
			\State{\textbf{Solve:} Unconstrained local LSQ problem}
		\EndIf
		\For{$j \gets 1$ to $\mathcal{\hat{N}}$}
		\Comment{where $\mathcal{\hat{N}}$ = nearest neighbors}
		\State{-- $\vec{P}_i(\Omega_i \cap (\Delta_j \cup \delta_j)) \rightarrow$ enqueue outgoing constraints}
		\State{-- Exchange shared DoFs with all nearest neighbor blocks}
		\State{-- $\vec{P}_i(\Delta_i \cup \delta_i)$ $\leftarrow$ dequeue incoming constraints}
		%
		%
		\State{-- Enforce constraints for $\vec{P}_i(\Delta_i \cup \delta_i) \cap \vec{P}_j, \forall j \in [1, \mathcal{N}]$:  Drive $\mathcal{F}_{ij} \rightarrow 0$}
		\EndFor
		
		\State{-- Update local error $E_i := \left\lVert \vec{Q}_i - R \vec{P}_i \right\rVert_{L_2} $}
		
		\EndFor

		\State{Check convergence metrics}
		
		\EndWhile
		\State{Write MFA to disk for analysis and visualization}
	\end{algorithmic}
\end{algorithm}

The DD techniques presented here for MFA computation are primarily implemented in Python, with main dependencies on \texttt{SciPy} for B-spline bases evaluations and linear algebra routines. Additionally, the drivers utilize Python bindings (\texttt{PyDIY}) for the \texttt{DIY}~\cite{morozov16} C++ library. \texttt{DIY} is a programming model and runtime
for block-parallel analytics on distributed-memory machines, built on \texttt{MPI-3}~\cite{dongarra13}.  Rather than programming
for process parallelism directly in \texttt{MPI}, the programming model in \texttt{DIY} is based on block parallelism. In \texttt{DIY}, data are decomposed
into subdomains called blocks. One or more of these blocks are assigned to processing elements (processes or threads) and the computation is
described over these blocks, and communication between blocks is defined by reusable patterns. 
\texttt{PyDIY} utilizes \texttt{PyBind11}~\cite{jakob17} and \texttt{MPI4Py}~\cite{dalcin11} to expose the interfaces in the C++ library. In our implementation, \texttt{PyDIY} is exclusively used to manage the data decomposition, including specifications to share an interface $\partial \Omega_{i,j}$ and ghost layers that represent the $\Delta \cup  \delta$ overlapping domains.

The overall approach is sketched in \algo{alg:pseudocode}.
We begin by decomposing the domain into a set of regular blocks aligned with the principal axes
of the global domain. Before enforcing constraints, the local subdomain solves are performed completely decoupled so that the discontinuous MFA to represent the partitioned input data is computed. 
The control point solution from this decoupled LSQ problem solver is then used as the DoF data that needs to be constrained with RAS iterative method.
We then begin iterating over the blocks to converge the shared DoFs through the linear constraints described in \sect{sec:solver-methodology}. 

At the start of each iteration, the control point constraints are exchanged between neighboring blocks in a regular nearest-neighbor communication pattern. This is sufficient to update the constraints $\vec{P}(\Delta \cup \delta)$. \texttt{DIY} 
sends and receives the constraint data to neighboring blocks based on the parallel data decomposition.
The nonlinear residual error in each subdomain is a function of the tensor product mesh resolution and degree $p$. At convergence, we expect to recover the subdomain error that is identical to the single subdomain case.

The final result, as described in \algo{alg:pseudocode}, is a global MFA that retains high-order continuity and accuracy of a single subdomain solve, but with excellent parallel efficiency to reduce total time to solution as the number of subdomains increases.

\section{Results}
\label{sec:results}

To demonstrate the effectiveness of the iterative algorithm for MFA computation, we devised a series of analytical closed form functionals and utilized real-world scientific datasets in both 2- and 3-dimensions obtained from high-fidelity simulations. All runs shown in this section were performed using the Python drivers written specifically for this work using the \texttt{DIY} domain decomposition infrastructure.

%

\subsection{\dimension{1} Results}\label{sec:results-1d}

In this section, detailed analysis on the convergence and accuracy of various MFA continuity recovery approaches are presented. 

\subsubsection{Comparison of Clamped vs Floating Boundary Knots}

To demonstrate the choice of using floating knots vs the low-order ($C^0$) continuous clamped knots at subdomain boundaries, we choose an analytical closed form reference solution of the form:

\begin{equation}
	F(x) = sinc(x) + sinc(2x-1) + sinc(3x+1.5), \forall x \in \Omega=[-4, 4]
	\label{eqn:1d-asymmetric-sinc}
\end{equation}

\begin{figure}[htbp]
	\centering
	\subfloat[Input analytical \dimension{1} solution profile\label{fig:1d-reference}]{%
		\includegraphics[width=0.44\textwidth]{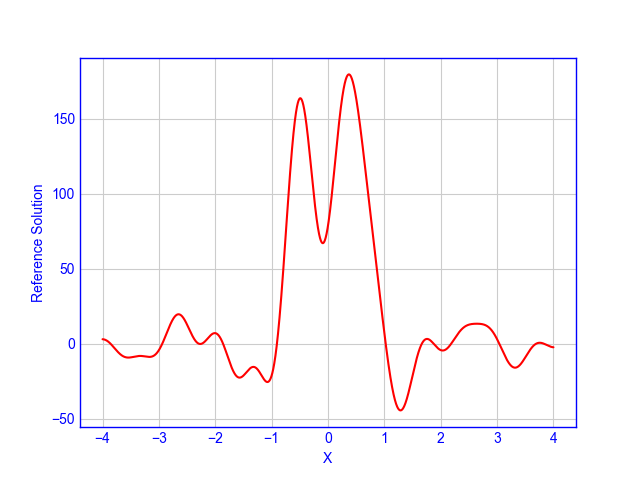}}
	\quad \quad
	\subfloat[Clamped $C^0$ continuous decoded solution\label{fig:1d-clamped}]{%
		\includegraphics[width=0.44\linewidth]{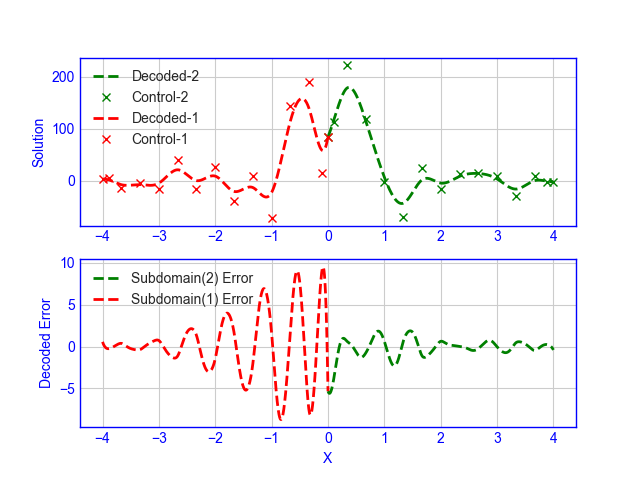}}
	\\
	\subfloat[Floating knots at $\Omega_{1,2}$ recovering $C^{p-1}$ continuity\label{fig:1d-unclamped}]{%
		\includegraphics[width=0.44\linewidth]{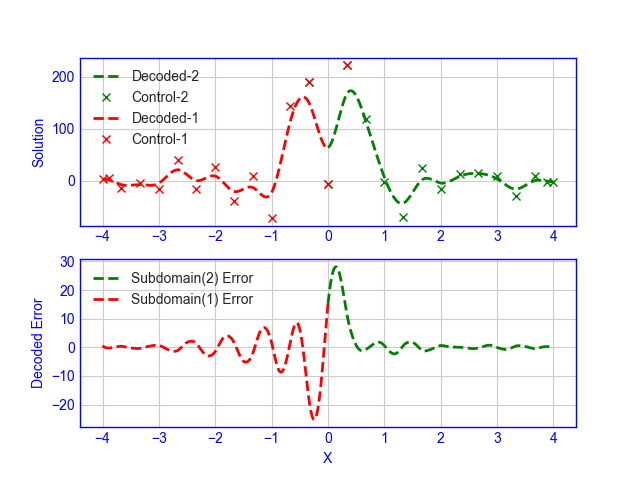}}
	\quad \quad
	\subfloat[Floating knots at $\Omega_{1,2}$ with $\left| \delta \right|=3$ \label{fig:1d-unclamped-aug}]{%
		\includegraphics[width=0.44\linewidth]{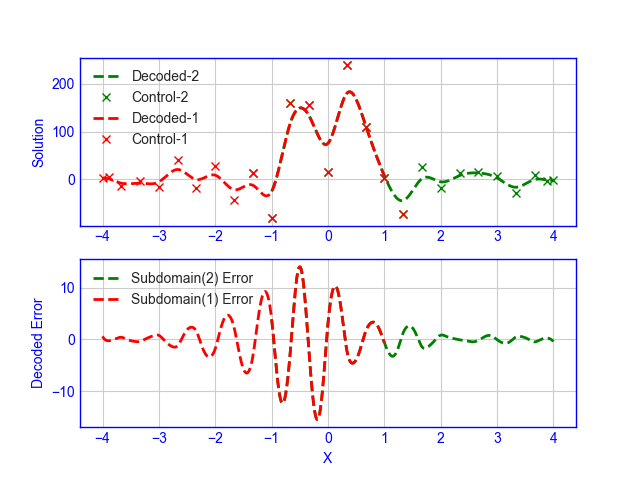}}
	\caption{\dimension{1} analytical sinc dataset with 10000 input points with $\mathcal{N}=2$ and $p=3$}
	\label{fig:comparison-clamped-floating}
\end{figure}

The reference solution $F(x)$, the results from the clamped knots, and floating knots with and without augmented overlap regions $\delta=p$ are shown in \fig{fig:comparison-clamped-floating}. The figures show the recovered solutions and the corresponding decoded error from MFA evaluation for a $\mathcal{N}=2$ and $p=3$ case. It is evident that the net error profile in the fully clamped subdivision in this example shows lower error as compared to the floating knot experiments. However, it is imperative to note that the former only shows $C^0$ regularity, while the floating knots fully recover high-order continuity at subdomain interfaces. Moreover, the use of augmented overlap regions ($\delta=3$) produce error profiles that resemble a single subdomain error profile in the domain, which is one of the key metrics of interest. These behaviors and conclusions extend to multi-dimensional setting as well.

\begin{figure}[htbp]
	\centering
	\includegraphics[width=0.44\textwidth]{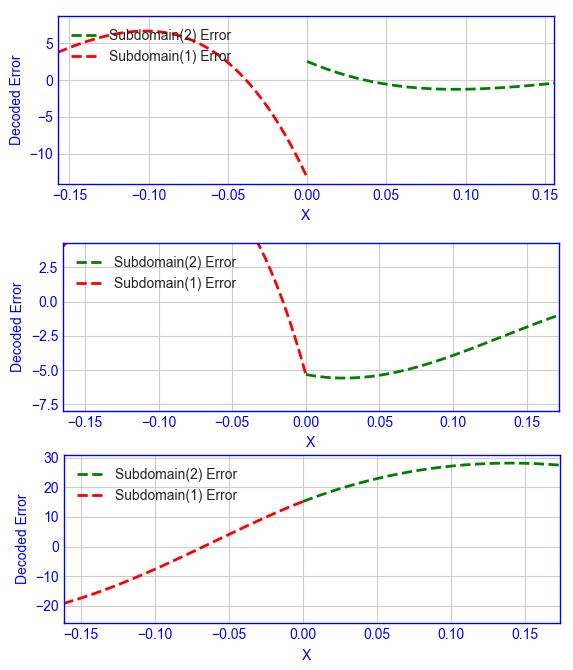}
	\caption{Zoomed error plots at interface $\Omega_{1,2}$ for \dimension{1} analytical dataset with $10^4$ input points, $\mathcal{N}=2$ and $p=3$.  Top: unconstrained and decoupled subdomain profile (discontinuous). Middle: clamped interface knots. Bottom: floating interface knots.}
	\label{fig:1d-error-comparison}
\end{figure}

To further demonstrate the continuity recovery behavior, we plot the error profile $\vec{E}$ for these approaches in \fig{fig:1d-error-comparison}, zoomed in around the interface $\Omega_{1,2}$. The unconstrained and decoupled LSQ solution procedure in the top shows that the reconstructed solution is discontinuous at the interface, as expected. Using the fully clamped approach that yields lower overall absolute error showcases only $C^0$ continuity at the interface, which may or not be sufficient depending on the use case utilizing the MFA representation. Finally, the bottom plot shows the smooth error profile from using the floating knots at the interface with full recovery of high-order continuity. We again emphasize that one could recover $C^0$ to $C^{p-1}$ continuity with this approach by choosing to use floating knots vs varying number of repeated knots at the interface.


%

\subsubsection{Error Convergence and Overlap Experiments}

\begin{figure}[htbp]
	\centering
	\subfloat[A double sinc profile in \dimension{1}\label{fig:error-1d-a}]{%
		\includegraphics[width=0.5\textwidth]{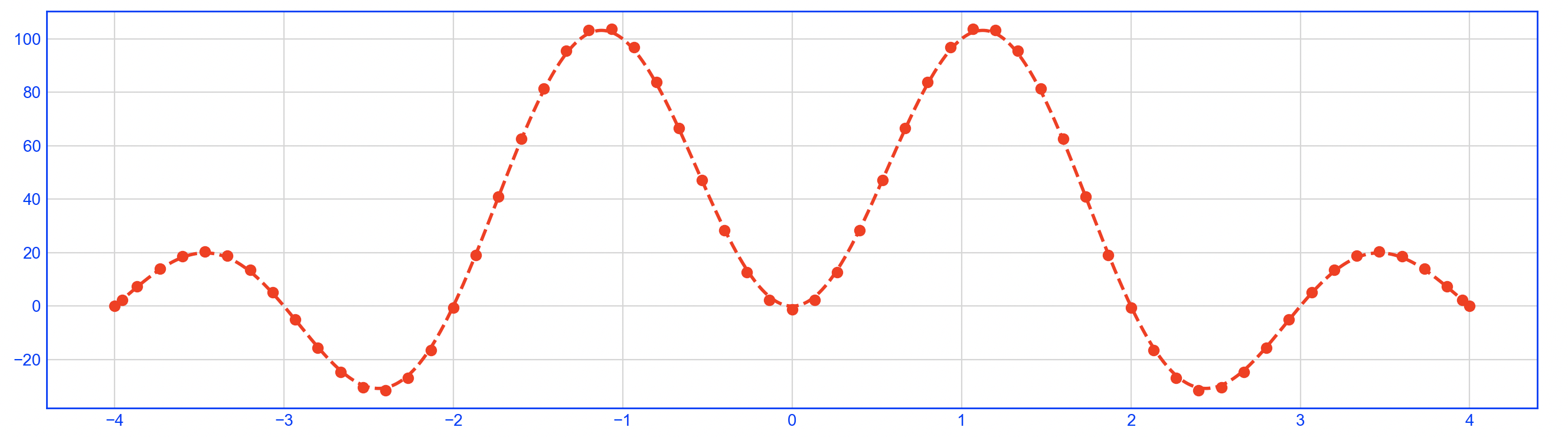}}
	\hfill
	\subfloat[Error profile for $\mathcal{N}=1$ and $p=3$\label{fig:error-1d-b}]{%
		\includegraphics[width=0.5\textwidth]{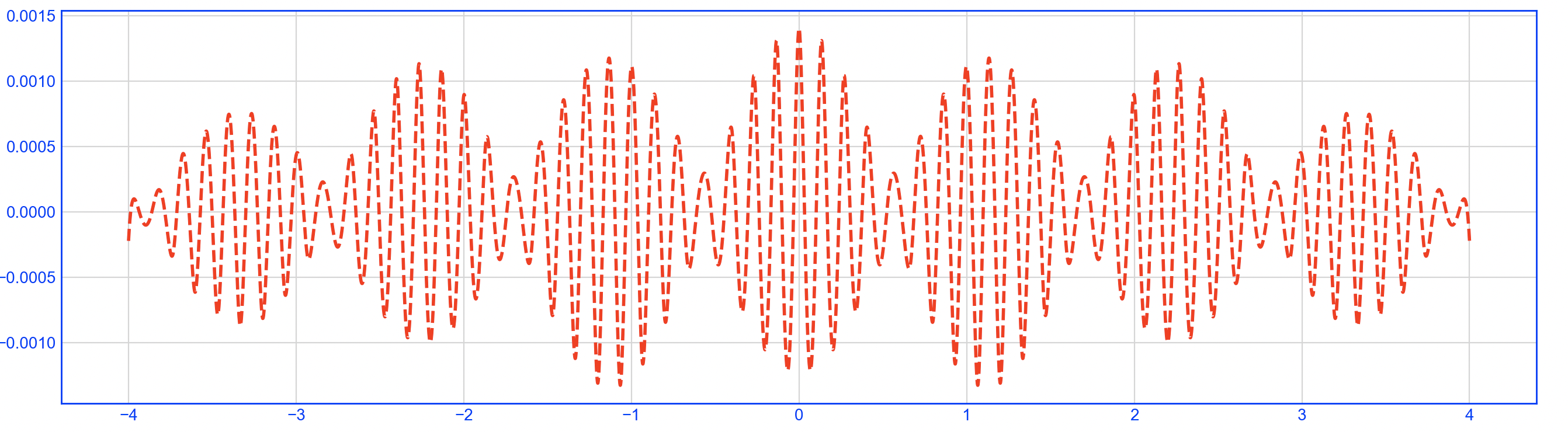}}
	\\
	\subfloat[Error profile for $\mathcal{N}=5$, $p=3$ and $\left| \delta \right|=0$\label{fig:error-1d-c}]{%
		\includegraphics[width=0.5\textwidth]{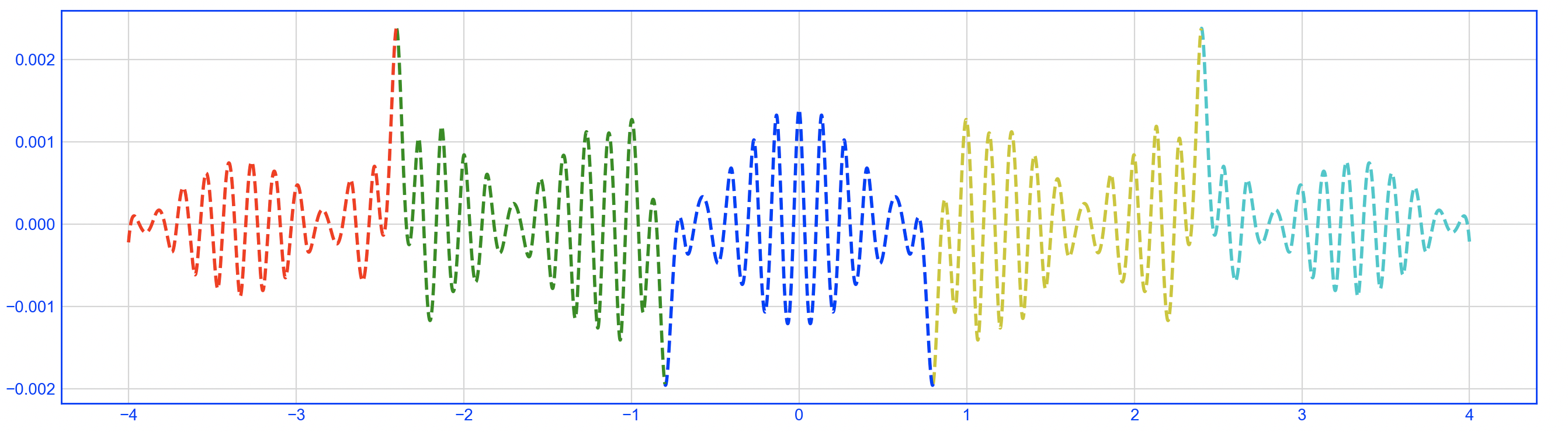}}
	\hfill
	\subfloat[Error profile for $\mathcal{N}=5$, $p=3$ and $\left| \delta \right|=p$\label{fig:error-1d-d}]{%
	\includegraphics[width=0.5\textwidth]{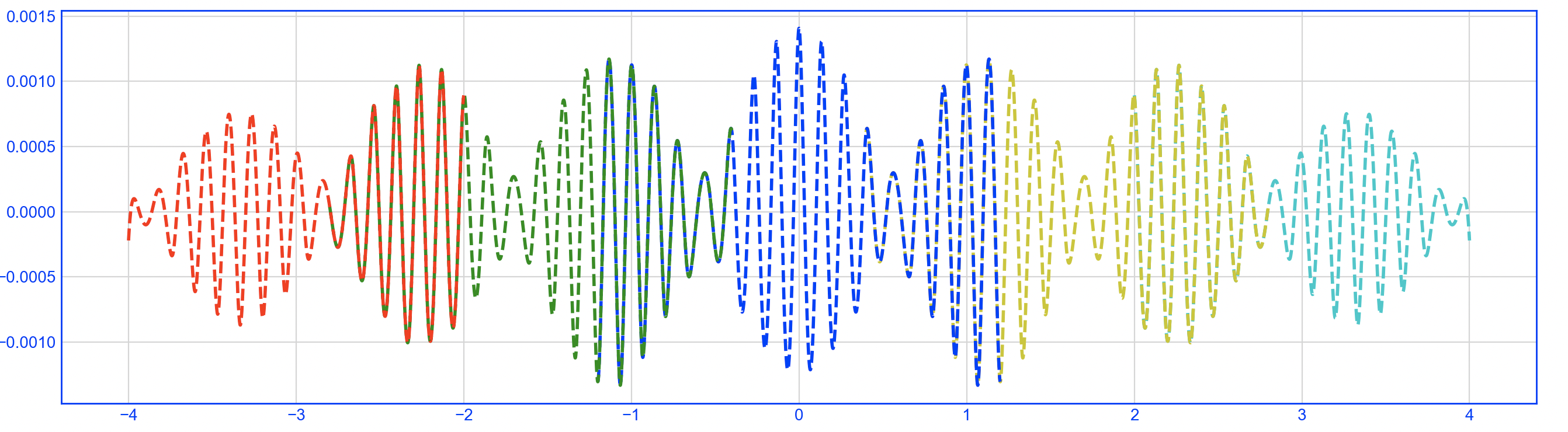}}
	\caption{Demonstration of error convergence, and effect of the overlapping spans to minimize numerical artifacts}
	\label{fig:error-1d}
\end{figure}

To determine the effect of using augmented or overlapped knot span regions ($\delta$) as the number of subdomains $\mathcal{N}$ are increased, we use a fully symmetric double-sinc function as shown below in \eqt{eqn:1d-symmetric-sinc} on a single subdomain as the reference solution as shown in \fig{fig:error-1d-a}, and with $\mathcal{N}=5$ for different values of augmented spans ($\left| \delta \right|=0$ and $\left| \delta \right|=3$).

\begin{equation}
	F(x) = sinc(x+1) + sinc(x-1), \forall x \in \Omega=[-4, 4]
	\label{eqn:1d-symmetric-sinc}
\end{equation}

It is evident from \fig{fig:error-1d-c} that when there is no augmented knot spans used in a multi-subdomain solver, the decoding of data at subdomain boundaries are influenced by contributions from both adjacent domain DoFs, which are enforced to be $C^{p-1}$ continuous by the constrained minimization solver. However, as we increase the number of overlap regions in terms of both the underlying data and the local bases support spans, the error profiles as shown in \fig{fig:error-1d-d} approaches the reference profile (with $\mathcal{N}=1$) shown in \fig{fig:error-1d-b}. Heuristically, for many of the problems tested, using $\left| \delta \right|=p$ provides optimal convergence as number of subdomains increase, even though increasing this parameter to $\left| \delta \right| =2p$ or higher will in general always improve the numerical accuracy at the cost of higher communication costs between neighboring subdomains.


\subsection{Multi-dimensional Problem Cases}

In this section, we present some results from applying parallel MFA for multidimensional problem cases.

\subsubsection{2D Problem}

\begin{figure}[htbp]
	\centering
	\subfloat[Nek5000 2D slice: Reference\label{fig:nek2d-reference-a}]{%
		\includegraphics[width=0.33\textwidth]{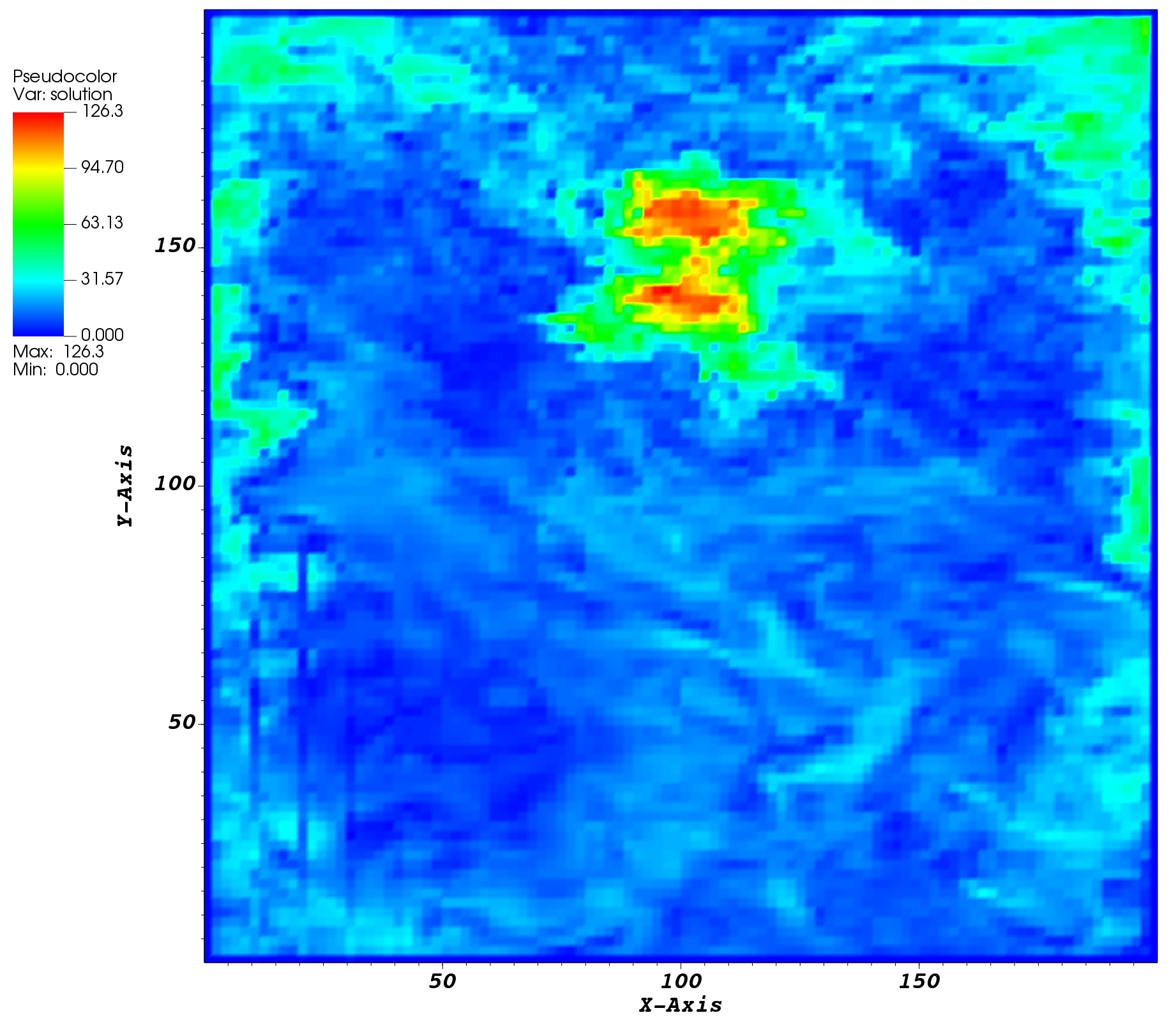}}
	\subfloat[Solution: 4x compression \label{fig:nek2d-d6-n5}]{%
		\includegraphics[width=0.33\textwidth]{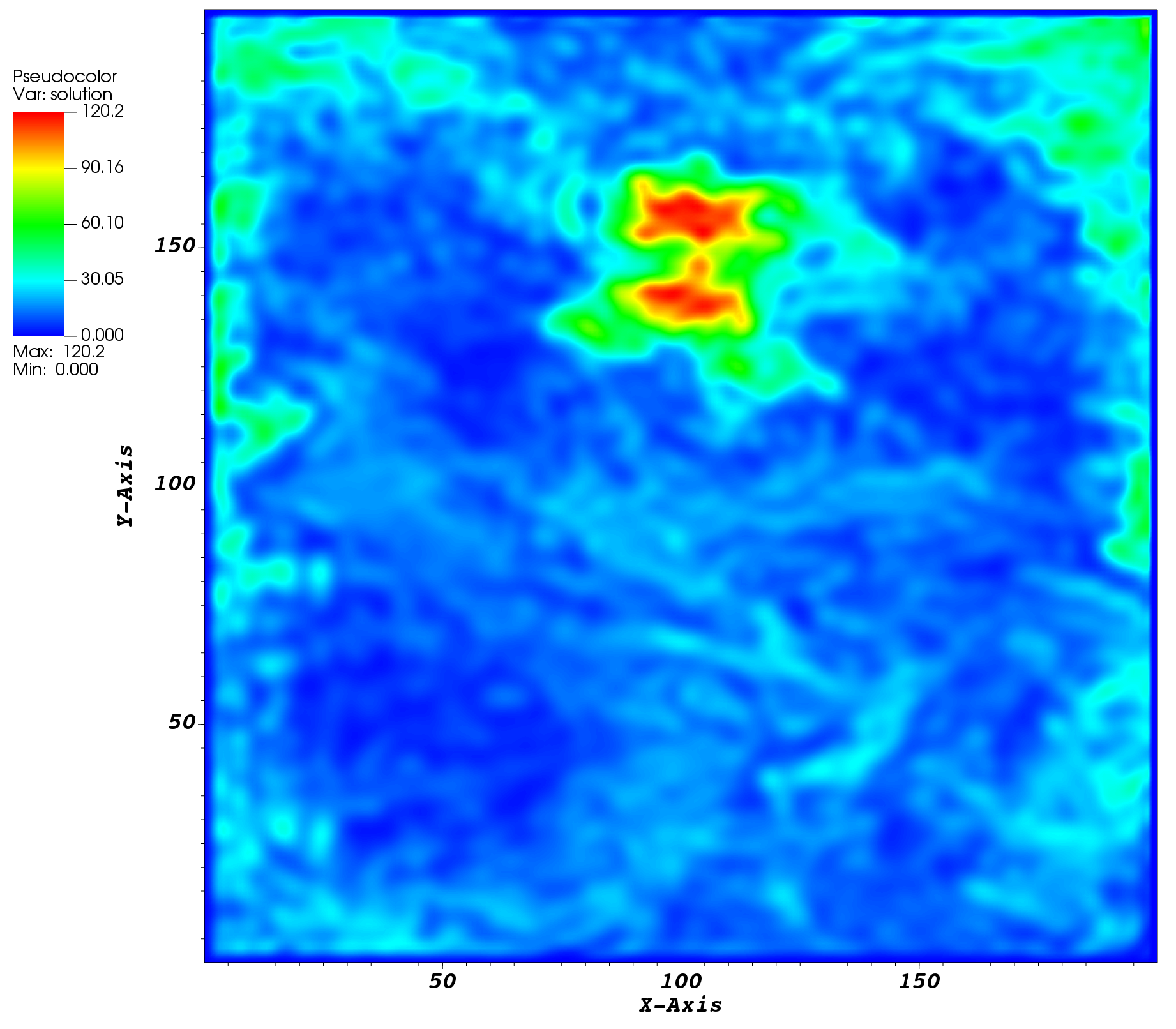}}
	\subfloat[Solution: 1x compression\label{fig:nek2d-d6-n1}]{%
		\includegraphics[width=0.33\textwidth]{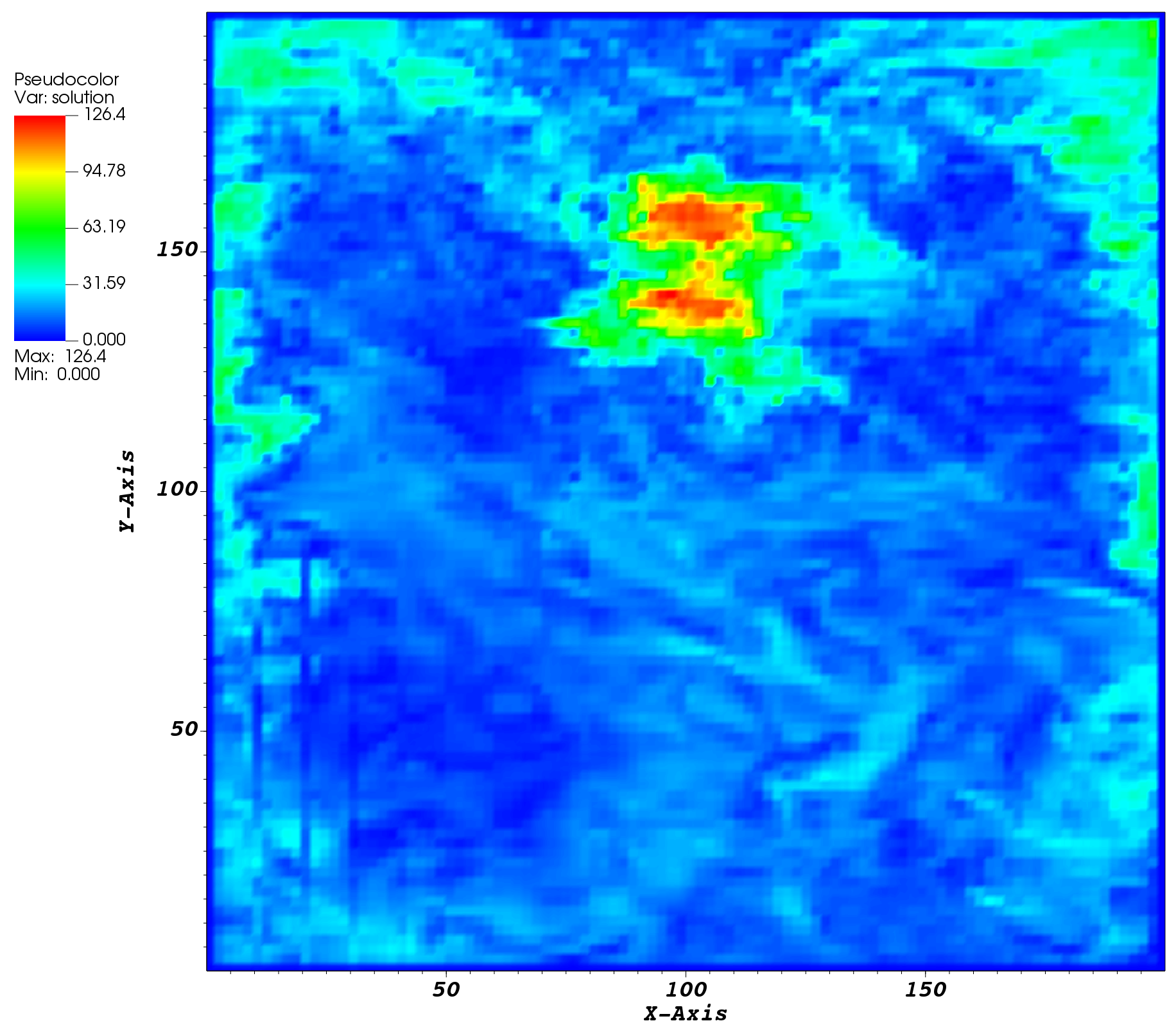}}
	\caption{2D slice of the Nek5000 \dimension{3} dataset ($200 \times 200$): reference profile and B-spline MFA with $p=6$, $\mathcal{N}=5\times5$, $\left| \delta \right|=p$ for 4X compression (middle), and 1X compression (right).}
	\label{fig:nek2d-ref-reconstruct}
\end{figure}

A \dimension{3} vector field representing the numerical results of a large-eddy simulation of Navier-Stokes equations for the MAX experiment \cite{merzari2010numerical} is representative of turbulent mixing and thermal striping that occur in the upper plenum of liquid sodium fast reactors. The data, generated by the Nek5000 solver \cite{deville2002}, have been resampled from their original topology onto a $200 \times 200 \times 200$ regular grid, and the magnitude of the velocity vector is associated with each \dimension{3} domain point \cite{peterka-mfa}.
Out of this dataset, a 2D slice (with $\left| \Omega \right|=200\times 200$) along the midplane in axial direction is used for our first study here. The reference solution and the converged, reconstructed solution with $\mathcal{N}=5 \times 5=25$ subdomains with $p=6$ and $\left| \delta \right|=2p$ is shown in \fig{fig:nek2d-ref-reconstruct} for different compression ratios. Depending on the usecase for MFA reconstruction, the converged error norms with 20 floating control points per subdomain yielding a net compression of 4X is sufficient to evaluate continuous derivatives everywhere in the domain $\Omega$. A full MFA representation (with compression ratio=1) is also shown in \fig{fig:nek2d-d6-n1}, which can fully reconstruct the original features in the input dataset in contrast to \fig{fig:nek2d-d6-n5} that shows a lossy smoothing of the sharp features in the original data.

\subsubsection{\dimension{3} Problem}

\begin{figure}[htbp]
	\centering
	\subfloat[S3D dataset profile\label{fig:s3d-refprofile}]{%
		\includegraphics[width=0.33\textwidth]{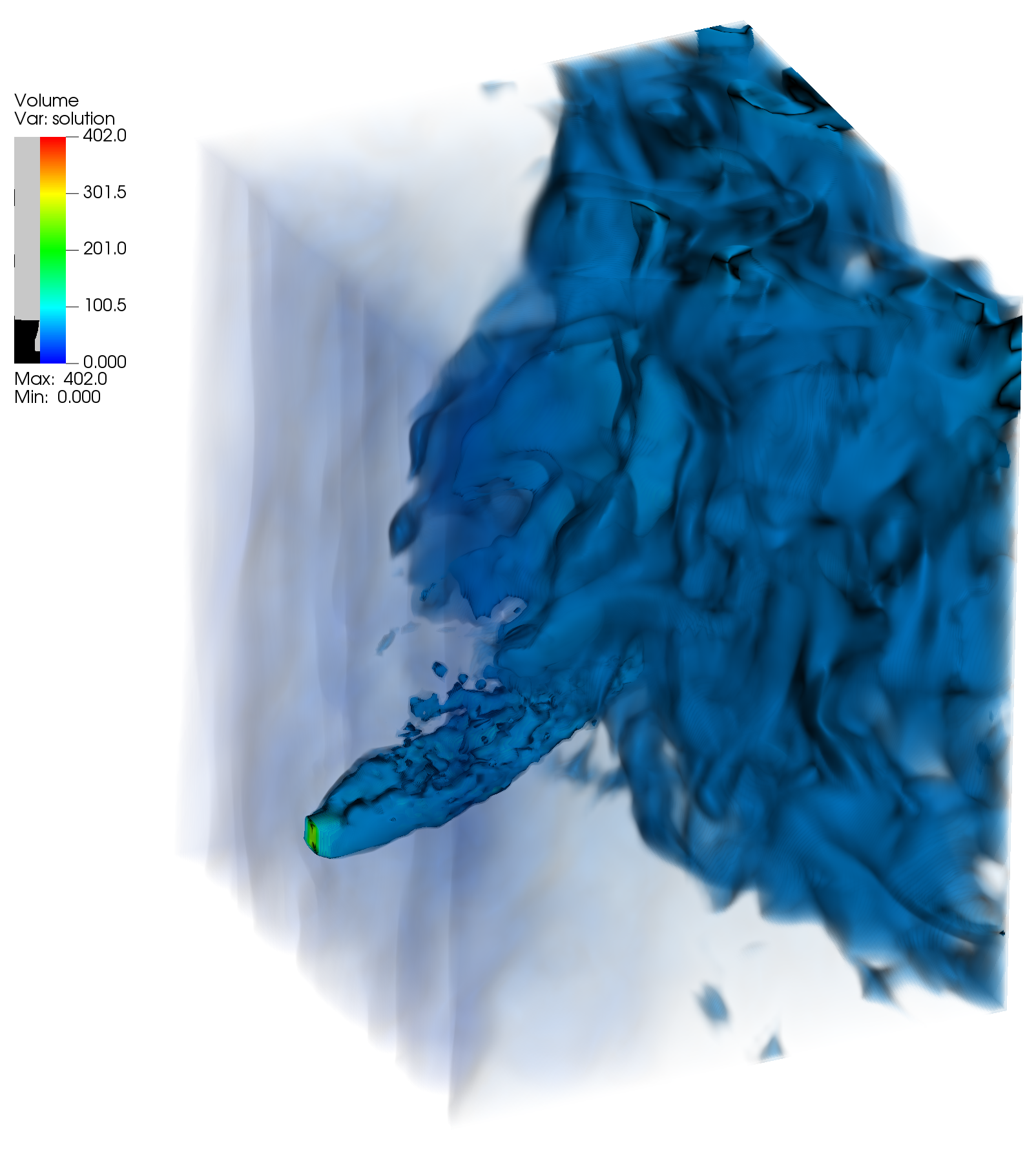}}
	\subfloat[Converged solution profile\label{fig:s3d-profile}]{%
			\includegraphics[width=0.33\textwidth]{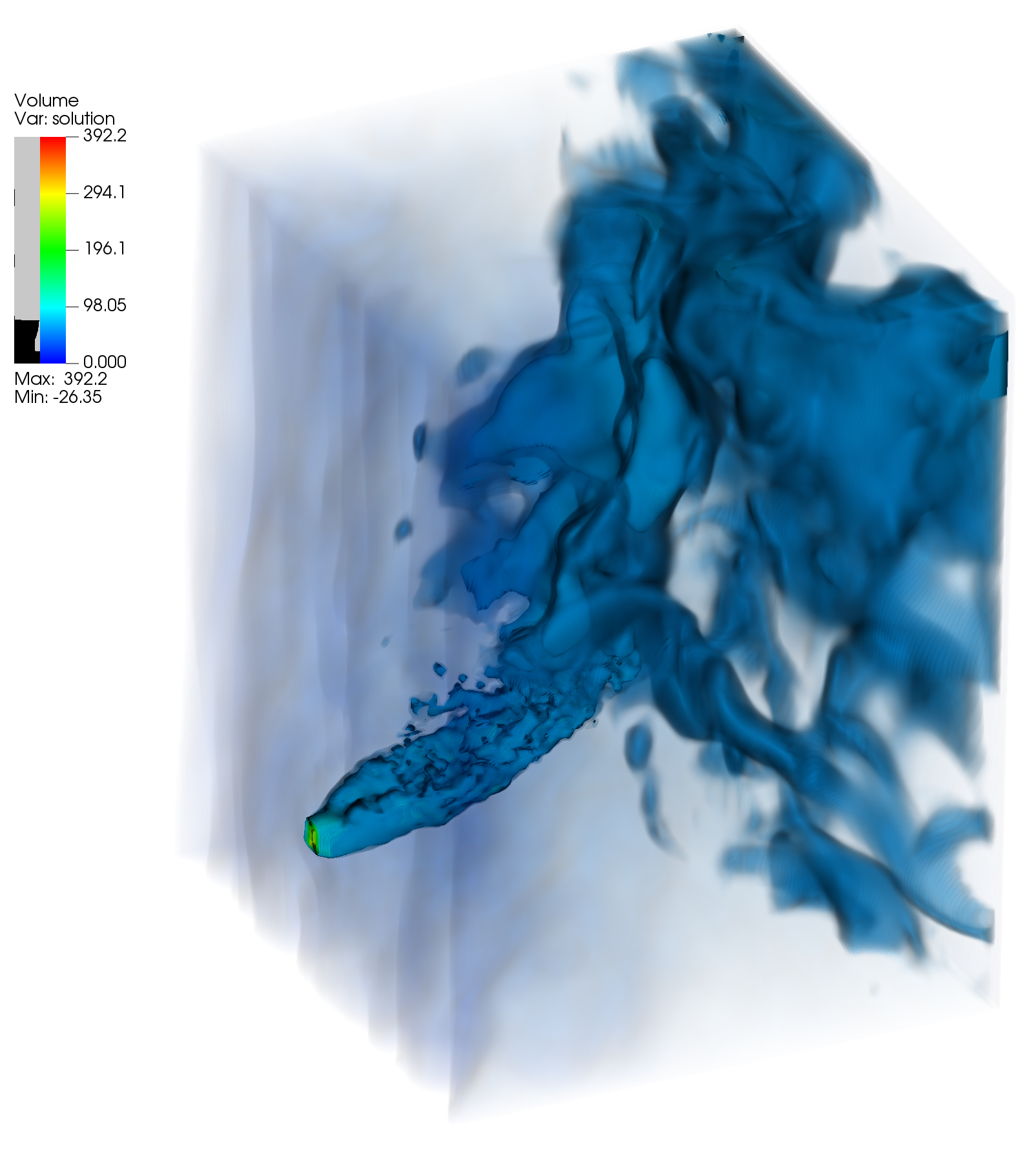}}
	\subfloat[Final error profile\label{fig:s3d-error}]{%
	\includegraphics[width=0.33\textwidth]{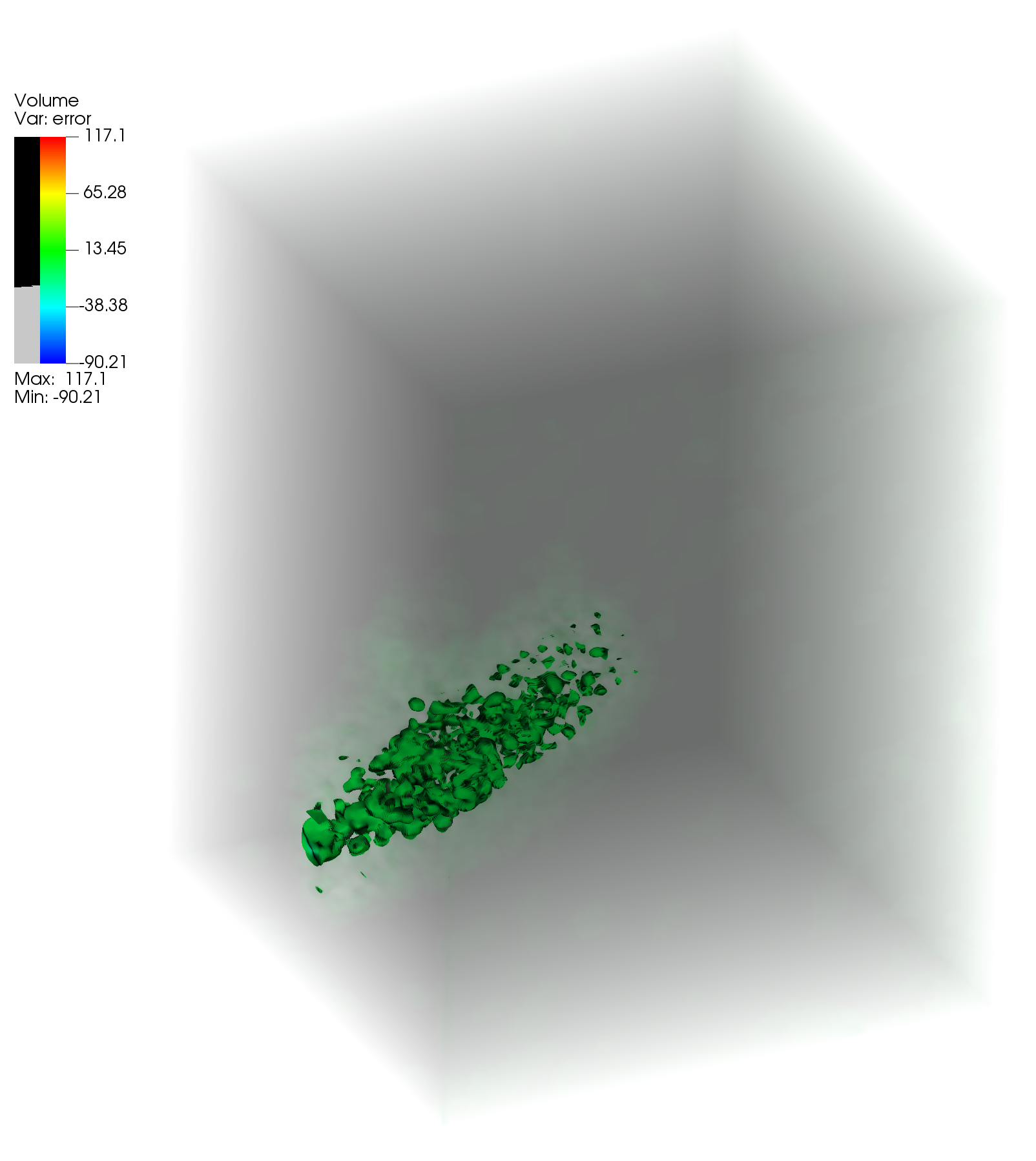}}
	\caption{Volume rendered S3D dataset with compression rate $\eta=9.5$: reference profile, converged MFA decoded profile and the corresponding reconstruction error with $\mathcal{N}=8\times 8\times 8=512$, $p=3$ and  $\left| \delta \right|=2p=6$}
	\label{fig:s3d-ref-error}
\end{figure}

Next, we present S3D, a turbulent fuel jet combustion dataset generated from a simulation in the presence of an external cross-flow \cite{chen2009}. The \dimension{3} domain has the span $\left| \Omega \right|=704 \times 540 \times 550$, with the raw data containing three components of the vector field. We choose to use the magnitude of this velocity field in our reconstruction study shown below in \fig{fig:s3d-ref-error} with 209M points. The converged MFA reconstruction shown in the figure with $8^3=512$ subdomains and $n=35$ per direction  in each subdomain yields a net compression ratio of $\frac{(704\times540\times550)}{(8\times35)^3} \approx 9.5$. While uniform refinement in knot spans does yield sufficient error reductions in most subdomains, utilizing adaptive error resolution with knot insertions and removals for MFA as previously used here for single subdomains \cite{nashed-rational} for the S3D problem can provide better reconstructions in addition to the iterative scheme introduced here. Such extensions will be pursued in the future. However, the experiments demonstrate that the MFA computations produce reconstructed data and  numerical errors that are consistent and convergent for arbitrary $\mathcal{N}, n, p$ and values of $\left| \delta \right|$, proving the feasibility of the algorithm.

\subsubsection{Error Convergence}


We utilize synthetic datasets shown in \eqt{eqt:2d-synthetic} and \eqt{eqt:3d-synthetic} to perform error convergence studies and to determine areas of maximal error that iteratively are resolved between neighboring subdomains.

%
%
%

\begin{eqnarray}
	F(x,y)  &=& sinc(\sqrt{x^2+y^2}) + sinc(2(x - 2)^2 + 2(y + 2)^2), \forall (x,y) \in [-4,4]^2, \label{eqt:2d-synthetic}\\
	\quad \quad F(x,y,z) &=& sinc(\sqrt{x^2 + y^2 + z^2}) + sinc(2(x - 2)^2 + (y + 2)^2 + (z - 2)^2), \forall  (x,y,z) \in [-4,4]^d.
	\label{eqt:3d-synthetic}
\end{eqnarray}

\begin{figure}[htbp]
	\centering
	\subfloat[\dimension{2} problem with $\mathcal{N}=5\times 5$: Iteration 1\label{fig:2d-iter1}]{%
		\includegraphics[width=0.44\textwidth]{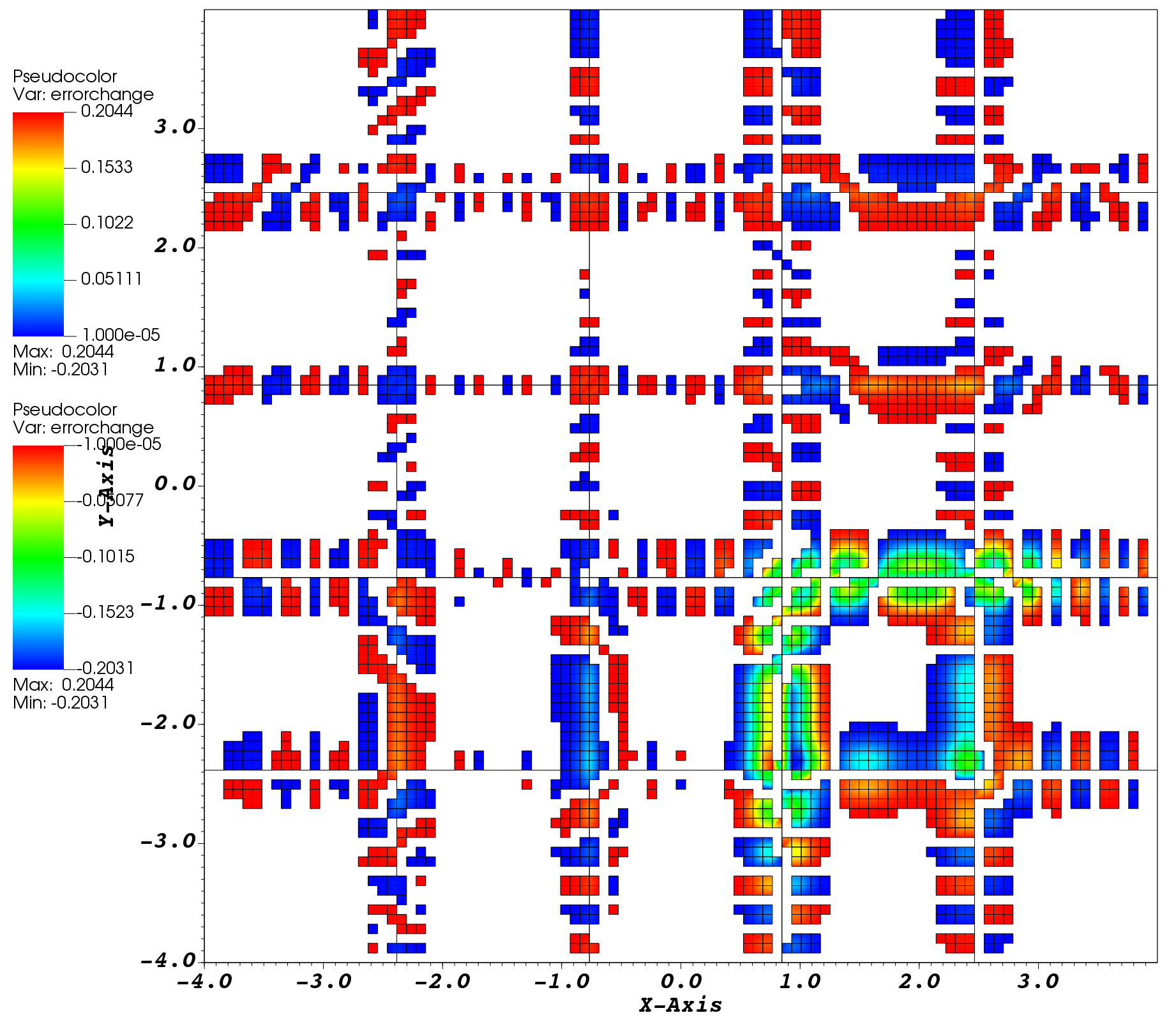}}
	\quad \quad
	\subfloat[\dimension{2} problem with $\mathcal{N}=5\times 5$: Iteration 2\label{fig:2d-iter2}]{%
		\includegraphics[width=0.44\textwidth]{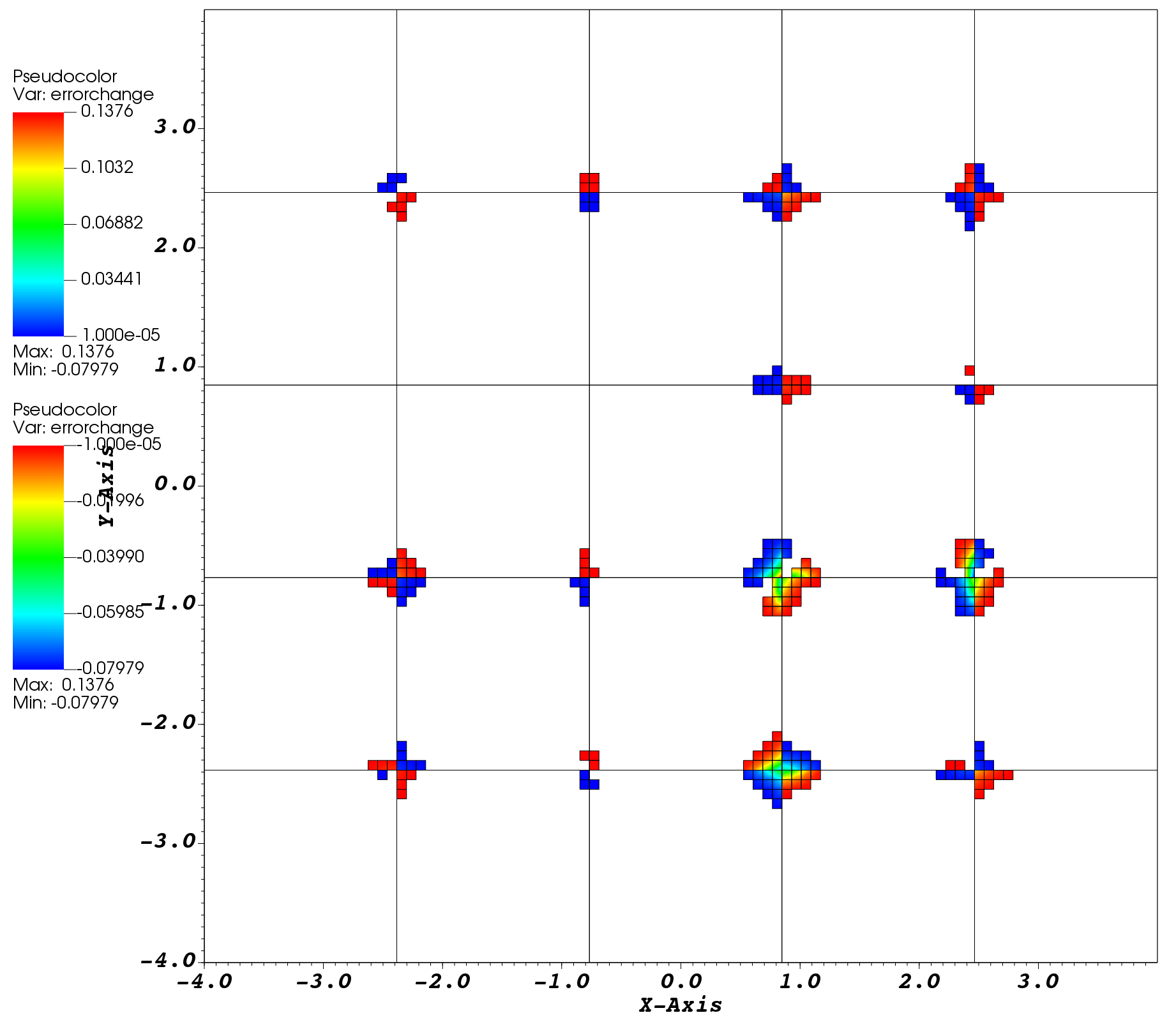}}\\
	\subfloat[\dimension{3} problem with $\mathcal{N}=3\times 3 \times 3$: Iteration 1\label{fig:3d-iter1}]{%
	\includegraphics[width=0.44\textwidth]{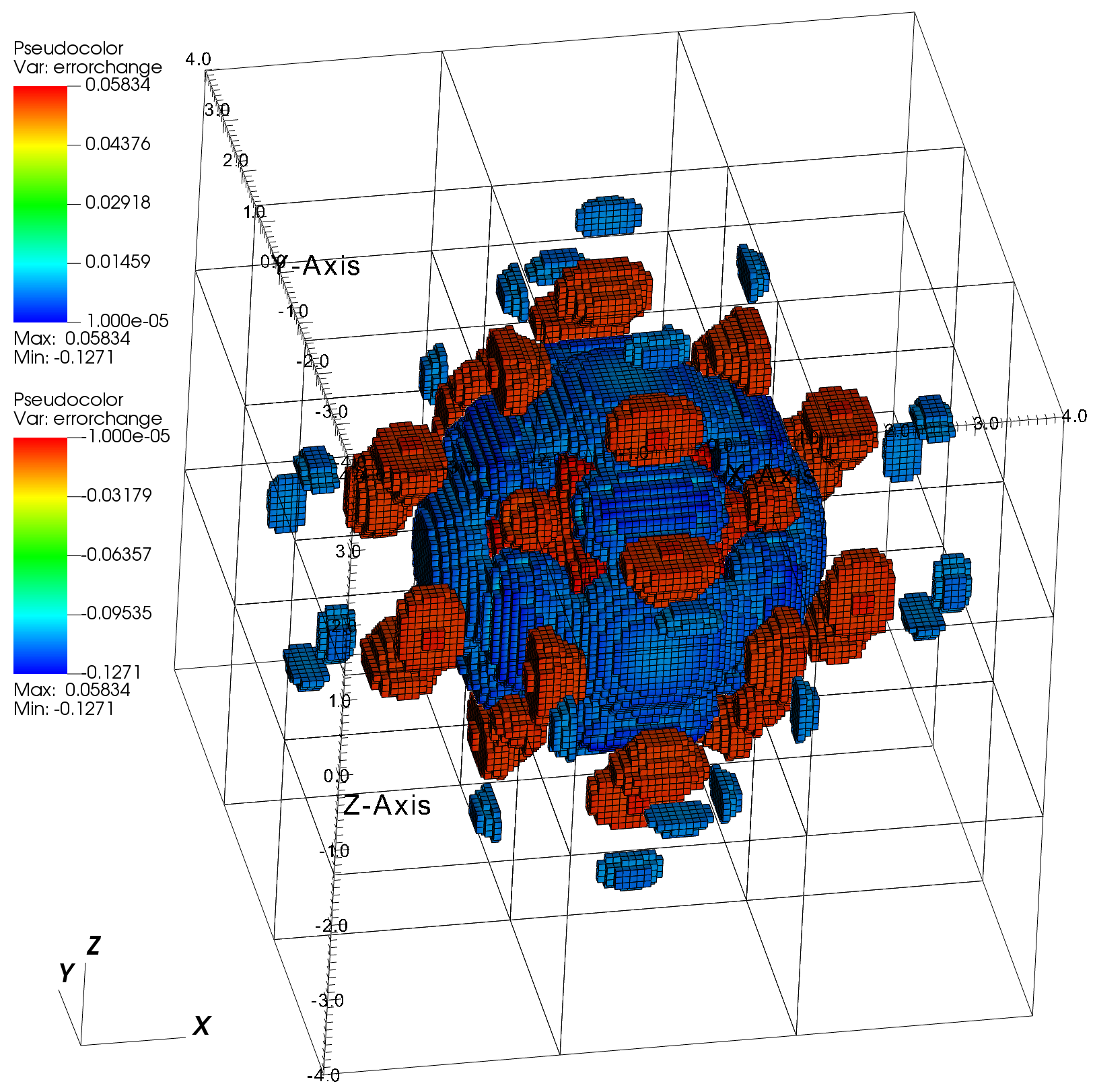}}
	\quad \quad
	\subfloat[\dimension{3} problem with $\mathcal{N}=3\times 3 \times 3$: Iteration 2\label{fig:3d-iter2}]{%
	\includegraphics[width=0.44\textwidth]{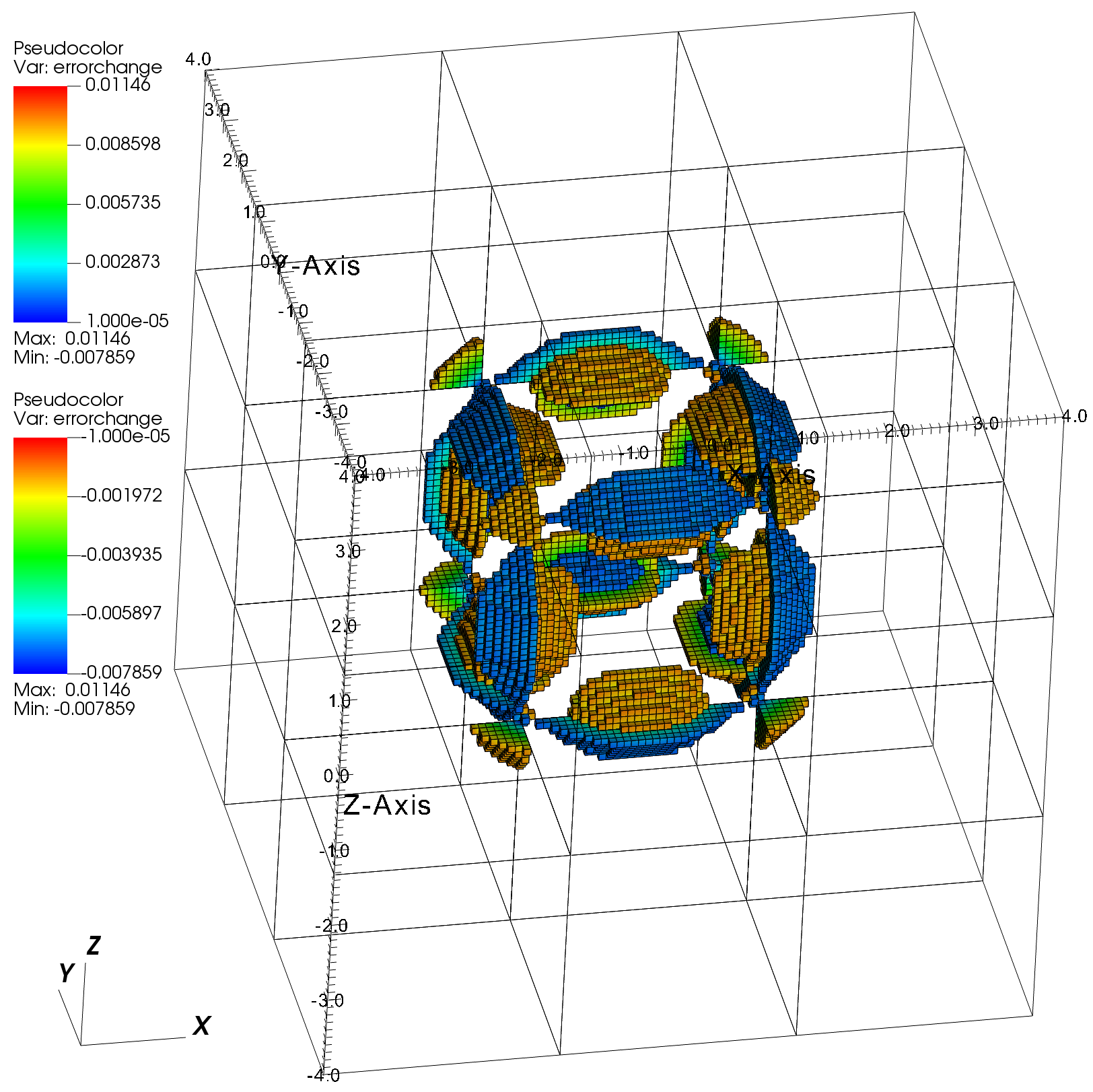}}\\
	\caption{Error convergence for \dimension{2} and \dimension{3} problems on the first and second iterate}
	\label{fig:2d-3d-error-convergence}
\end{figure}

We plot the change in error $\vec{E}$ between subsequent iterations of the RAS scheme for both \dimension{2} and \dimension{3} problem cases with $\left| \delta \right|=0$ in \fig{fig:2d-3d-error-convergence}. This clearly demonstrates that the $\mathcal{SS}$ interface values are resolved first, and then $\mathcal{MS}$ DoFs are resolved further. In all cases, the iterations converge in 2 steps, independent of $\mathcal{N}$ or $\delta$.

\subsubsection{Parallel Scalability}\label{sec:parallel-scalability}

To demonstrate the parallel performance of the implemented RAS iterative scheme for MFA computation with continuity preservation, we employed both closed-form synthetic and real simulation datasets in \dimension{2} and \dimension{3}. In the following sections, we present both strong scaling and weak scaling studies performed on the Theta Cray XC40 supercomputer operated by the Argonne Leadership Computing Facility (ALCF), which provides 4,392 KNL compute nodes with 64 compute cores and 192 GB DDR4 RAM per node. The interconnect is based on the Aries Dragonfly high speed network.

\subsubsection*{Strong Scaling}

We consider both \dimension{2} and \dimension{3} problems to demonstrate the scaling behavior of the presented parallel MFA computational algorithm in \algo{alg:pseudocode}. One key consideration that drove selection of the subdomain sizes, and the floating knot span descriptions, is motivated by the metric to recover the original error profile from a single subdomain case. 
Verification studies were performed during this strong scaling test to ensure that the local subdomain errors computed on a single task, and on different process counts remain the same at convergence. This verification is important to reiterate the fact that the approximation error due to the constrained solves to recover higher-order continuity does not significantly affect the error metrics for the MFA as $\mathcal{N}$ increases. For this reason, we used synthetic datasets generated with closed form equations for \dimension{2} and \dimension{3} studies as shown in \eqt{eqt:2d-synthetic} and \eqt{eqt:3d-synthetic}.



\begin{figure}[htbp]
	\centering
	\subfloat[\dimension{2} strong scaling\label{fig:strong-scaling-2d-a}]{%
		\includegraphics[width=0.49\textwidth]{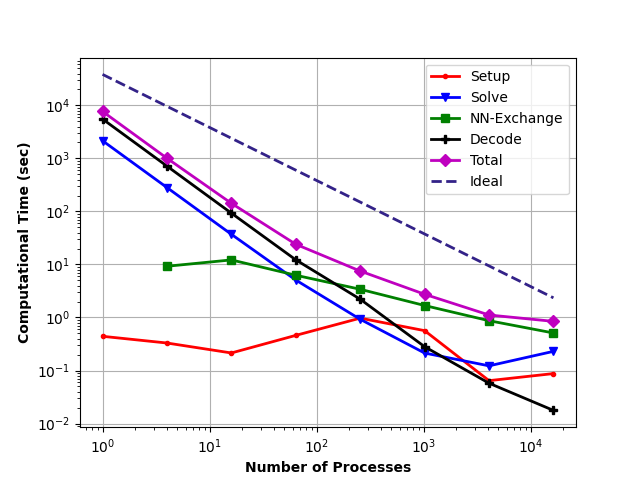}}
	\subfloat[\dimension{2} strong scaling with $\left| \delta \right|=p$\label{fig:strong-scaling-2d-b}]{%
		\includegraphics[width=0.49\textwidth]{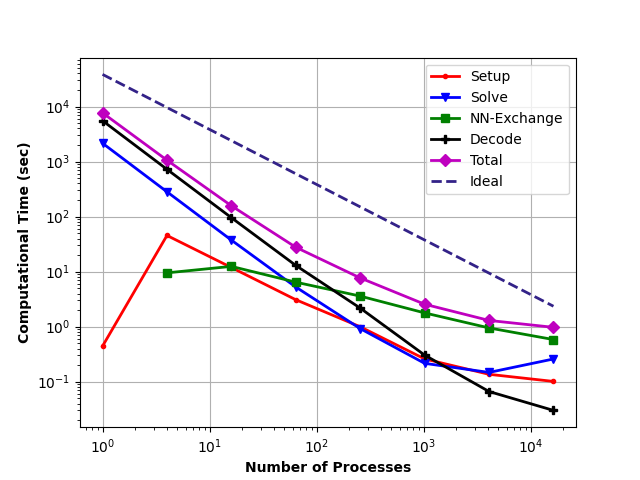}}\\
	\subfloat[\dimension{3} strong scaling\label{fig:strong-scaling-3d-a}]{%
		\includegraphics[width=0.49\textwidth]{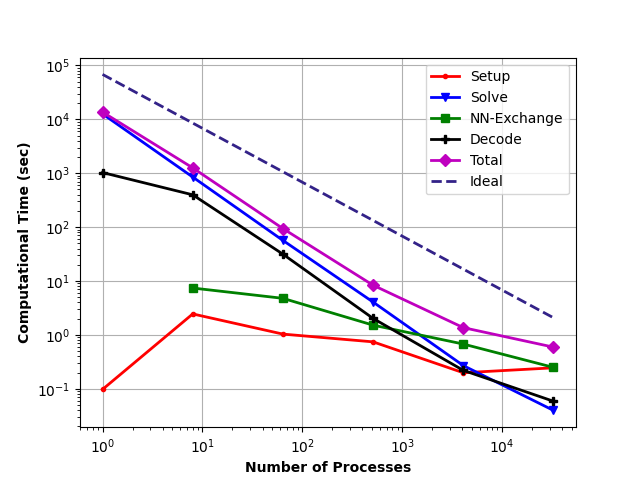}}
	\subfloat[\dimension{3} strong scaling with $\left| \delta \right|=p$\label{fig:strong-scaling-3d-b}]{%
		\includegraphics[width=0.49\textwidth]{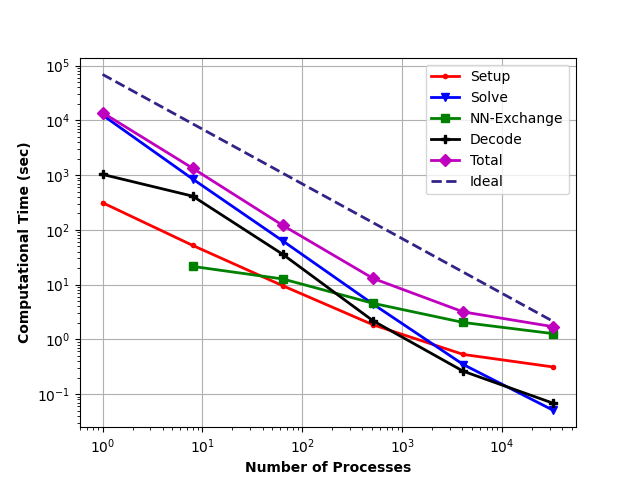}}\\
	\caption{Strong scaling performance of RAS solver with $p=3$}
	\label{fig:strong-scaling}
\end{figure}

The strong scaling tests were performed on 1 to 16,384 parallel tasks in \dimension{2}, increasing by a factor of $2^2=4$, and the \dimension{3} tests were executed on 1 to 32,768 tasks, increasing by a factor of $2^3=8$. We also note that the \dimension{2} and \dimension{3} case for these studies used 400M and 1.331B input points respectively, with a corresponding $\eta=4$ and $\eta=1.25$ that was maintained constant for all runs. 
In order to also better understand the effects of using augmented overlap regions ($\delta$) on scalability, two cases with the choice of $\left| \delta \right|=0$ and $\left| \delta \right|=p$ are shown in \fig{fig:strong-scaling}. 
%
%
The Python driver utilized \texttt{DIY} to handle block decompositions and rank assignments, as the total number of tasks used in the parallel run was increased. We measured the overall computational time for setting up the problem, the initial subdomain solves, and the consequent RAS iteration cycle to convergence, which includes the nearest neighbor communication at each iteration. We also show the time for decoding the MFA that is used to measure the errors in each subdomain, and the overall total that includes the effort spent on each of these various tasks. This task-wise breakdown helps us clearly visualize the steps that scale linearly and the ones that do not.



As expected, the RAS iterative scheme shows excellent scalability for the chosen datasets, and the overall time to compute the MFA in parallel was reduced at a nearly ideal rate up to $10^3$ MPI tasks as $\mathcal{N}$ increases, while ensuring $C^{p-1}$ continuity in the subdomain interfaces.  It is important to note that the dominant computational time is driven by the decoupled LSQ solution computation and decoding operations, which are embarassingly parallel as the size of the subdomains (determined by number of knot spans) decreases in direct proportion to the tasks. Given that the scalability of the linear algebraic LSQ solvers \cite{benzi2002, dalcin2011} and Sparse Matrix Vector (SpMV) products used in the decode tasks are well understood, the bottlenecks potentially occur only from the nearest neighbor communication for constraint data exchanges, which remain relatively small in magnitude up to 32K processes tested in the \dimension{3} experiments. The overall strong scaling efficiency remains around 50\% for both overlapping and non-overlapping \dimension{2} problem cases at 16,384 tasks. However, the added setup cost and less than ideally scaling nearest neighbor communications reduce the \dimension{3} problem efficiency for the overlapping subdomain cases to 25\% at 32,768 tasks from 70\% in the nonoverlapping cases. 

Additionally, the parallel efficiency degradation behavior in augmented runs due to the high setup cost occurs due to the current choice of implementations to determine intervals in both the input space data ($\vec{Q}$) and extra knot spans ($\delta$) that need to be exchanged with neighboring subdomains. This setup cost is purely local in the case when $\left| \delta \right|=0$, in contrast to the communication dominated setup when $\left| \delta \right| > 0$. 
It is evident that the cost of nearest neighbor exchanges at scale is driven by the size of the messages being transferred between subdomains, which is a function of $\left| \delta \right|$. The increased amount of data shows smaller degradation in communication for \dimension{2} problems as compared to \dimension{3}, and hence the strong scaling behavior is dependent on the dimension of the problem being solved that is representative of performance in all structured grid solvers for PDEs.

\subsubsection*{Weak Scaling}

\begin{figure}[htbp]
	\centering
	\subfloat[\dimension{2} weak scaling\label{fig:weak-scaling-2d}]{%
		\includegraphics[width=0.5\textwidth]{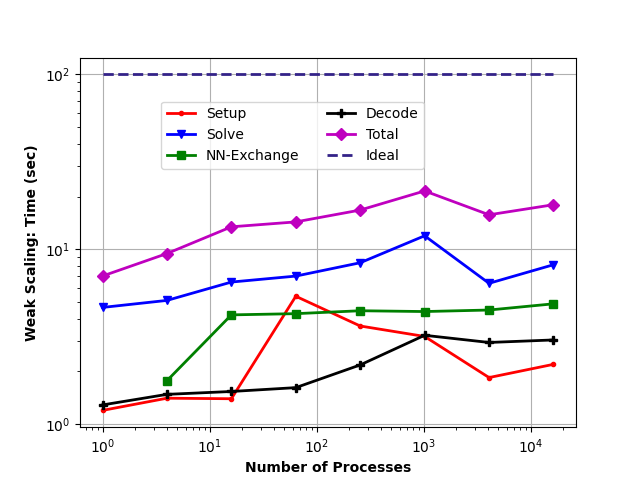}}
	\hfill
	\subfloat[\dimension{3} weak scaling\label{fig:weak-scaling-3d}]{%
		\includegraphics[width=0.5\textwidth]{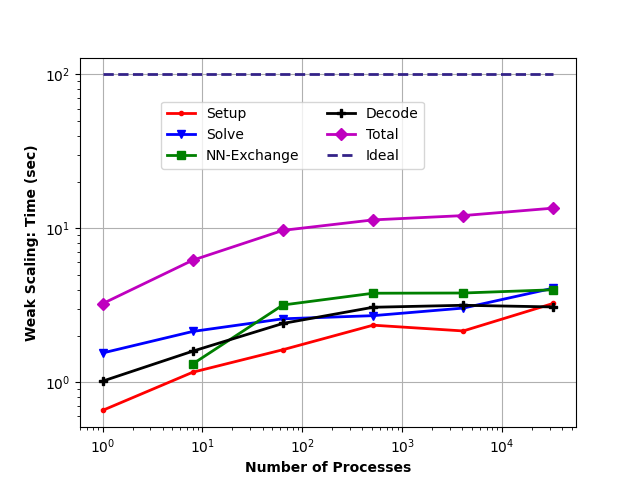}}\\
	\caption{Weak scaling performance of RAS solver with $p=3$ and $\left| \delta \right|=p$}
	\label{fig:weak-scaling}
\end{figure}

Given that the performance of the overlapping and augmented MFA scheme was comparable to non-overlapping cases ($\left| \delta \right|=0$), and since the error reduction from having extra overlaps always results in better solution recovery, we strictly focus on overlapping cases alone for the weak scaling study. Here, the overall work per subdomain is maintained constant, and the number of tasks are increased from 1 to 16,384 in \dimension{2}, and from 1 to 32,768 in \dimension{3}, similar to the strong scaling study.
The weak scaling results maintaining an overall MFA compression rate of $\eta=2^d,\forall d \in [2,3]$ is shown in \fig{fig:weak-scaling}.

The weak scaling study demonstrates that the overall performance of the RAS iterative scheme for large number of subdomains does not significantly affect the parallel efficiency, which are around 40\% in \dimension{2} and 24\% in \dimension{3} at the fine limit tested. The subdomain solve and nearest neighbor data exchange dominate the overall time to solution. However, it is important to note that the actual runtime for the MFA computation only grows by a factor of 2, even on 16K processes or more.

\subsubsection*{Performance Study on S3D Dataset}

\begin{figure}[htbp]
	\centering
	\includegraphics[width=0.6\textwidth]{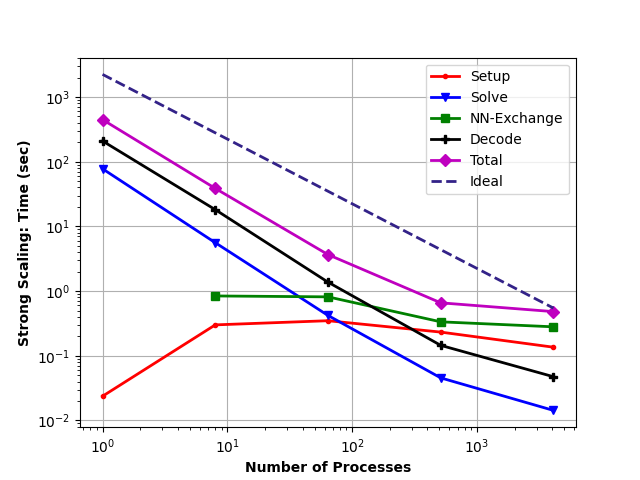}
	\caption{Strong scaling for the S3D dataset with $p=5$}
	\label{fig:s3d-strong-scaling}
\end{figure}

Finally, we consider the case of the S3D combustion dataset shown in \fig{fig:s3d-ref-error} and measure the strong scaling performance on up to 4,096 processes. Using parallel MPI-IO implemented with \texttt{DIY}, and exposed through the Python interface, a strong scaling performance study was measured on this realistic dataset and shown in \fig{fig:s3d-strong-scaling}. Note that the IO cost for reading the chunk of data required on each task is included in the setup time shown in the figure.

The performance and error analysis indicate good speedup to reduce overall time for MFA computation, until nearest neighbor communication and data exchanges start dominating the overall workflow. These results show similar behavior to the strong scaling studies performed on synthetic datasets and provide confirmation on the feasibility of the presented approach for tackling real-world large datasets.

\section{Summary}
\label{sec:conclusions}


We have presented a scalable DD approach to tackle the issue of discontinuous B-spline based MFA representations when performing the computations in parallel. The Restricted Additive Schwarz (RAS) method is a natural algorithmic fit for data analysis problems to create efficient MFA solutions in parallel. Through the use of Schwarz-based iterative schemes, combined with constrained local subdomain solvers, the two-level iterative technique has been shown to be robust in converging to the compressed functional representation of the given data, without sacrificing the approximation accuracy measured on a single subdomain of equivalent control point resolution. Replacing B-spline bases with NURBS bases ($W \ne 1$) only requires imposing the constraints on the $\vec{P}_i W_i$ data instead of $\vec{P}_i$ alone, which is naturally accomplished with minor modifications in \algo{alg:pseudocode}. This can also be combined with aposteriori error measures \cite{nashed-rational} to adaptively resolve solution variations while ensuring higher-order continuity across subdomain boundaries with appropriate knot insertion, removal and communication of shared DoFs in $\Delta \cup \delta$ regions. Another natural way to ensure continuity across adaptively resolved NURBS or B-spline patches would be to use T-splines \cite{sederberg-2004}, which are specifically designed for merging higher-dimensional surfaces with non-matching knot locations. All presented ideas should extend for T-splines instead of B-splines as well with modifications to $\mathcal{C}$  in \eqt{eq:global-constrained-problem}.

We have demonstrated that the use of overlap layers $\delta$ can certainly improve the overall MFA accuracy, with a slightly larger one-time setup cost that gets amortized in the overall computation time. We determined that for all the problems tested, including real datasets, $\left| \delta \right|=p$ to $\left| \delta \right|=2p$ is optimal in terms of error recovery and computational cost even for \dimension{3} problems up to 32,768 tasks.
%
%
The iterative scheme shows good parallel performance for both \dimension{2} and \dimension{3} problems tested, and the parallel efficiency degrades only when the cost of nearest neighbor subdomain data exchanges start to creep up beyond the cost of the local constrained subdomain solve. Given that scaling characteristics of these processes are well understood in the literature, the parallel speedups behave predictably well at scale on large computing machines tested.

The \texttt{PyDIY} based Python implementations for 1-, 2- and 3-dimensional problems have been shown to resolve large, complex solution profiles with strong gradient variations, even under decreasing subdomain sizes. Depending on the needs for visualization or in-situ analysis, choices on clamped or floating knots can be made with no modifications to the implementation. This scheme can also be used to achieve scalable high-order solution field transfers between component models, a process more generally referred to as {\em remapping} \cite{dukowicz1987accurate}, by imposing constraints on the subdomain solvers to satisfy various metrics of interest \cite{mahadevan2022} such as global conservation and monotonicity. The exploration of parallel MFA for such applications will be pursued in the future.

\section*{Acknowledgments}

This work is supported by Advanced Scientific Computing Research, Office of Science, U.S. Department of Energy, under Contract DE-AC02-06CH11357, program manager Laura Biven. We gratefully acknowledge the computing resources provided on the Theta supercomputer operated by Argonne Leadership Computing Facility, which is a DOE Office of Science User Facility supported under Contract DE-AC02-06CH11357, which was used to generate all the strong and weak scaling results. This research also used resources in the Bebop, a high-performance computing cluster operated by the Laboratory Computing Resource Center (LCRC) at Argonne National Laboratory to test and optimize performance of the drivers at scale.

%

\bibliographystyle{siamplain}
\bibliography{asm-mfa}

\pagebreak
\section*{Government License}

The submitted manuscript has been created by UChicago Argonne, LLC, 
Operator of Argonne National Laboratory (“Argonne”). Argonne, a U.S. 
Department of Energy Office of Science laboratory, is operated under 
Contract No. DE-AC02-06CH11357. The U.S. Government retains for itself, and others acting on its behalf, a paid-up nonexclusive, irrevocable 
worldwide license in said article to reproduce, prepare derivative 
works, distribute copies to the public, and perform publicly and display 
publicly, by or on behalf of the Government. The Department of Energy 
will provide public access to these results of federally sponsored 
research in accordance with the \href{http://energy.gov/downloads/doe-public-access-plan}{DOE Public Access Plan}.

\end{document}